\documentclass[11pt,reqno,a4paper,]{article}

\usepackage[utf8]{inputenc}
\usepackage[T1]{fontenc}
\usepackage{lmodern}

\usepackage{float}
\usepackage{yhmath}
\usepackage{centernot}
\usepackage{amsmath,amssymb,amsthm,mathtools}
\usepackage{bm}
\usepackage{enumitem}
\usepackage{graphicx}
\graphicspath{{figures/}}
\usepackage{caption}
\usepackage{subcaption}
\usepackage[colorlinks=true,linkcolor=blue,citecolor=blue,urlcolor=blue]{hyperref}
\usepackage{algorithm}
\usepackage{algpseudocode}
\usepackage{booktabs}
\usepackage{geometry}
\newtheorem{proposition}{Proposition}[section]

\newtheorem{theorem}{Theorem}[section]
\newtheorem{lemma}{Lemma}[section]
\newtheorem{corollary}{Corollary}[section]
\theoremstyle{definition}
\newtheorem{definition}{Definition}[section]
\newtheorem{remark}{Remark}[section]

\renewcommand{\thefootnote}{\fnsymbol{footnote}}
\newcommand{\precprec}{\prec\mathrel{\mkern-5mu}\prec}
\newcommand{\IR}{\mathbb{IR}}
\newcommand{\R}{\mathbb{R}}

\algrenewcommand\alglinenumber[1]{\textbf{Step #1.}}

\title{\textbf{Proximal gradient method for multiobjective interval valued optimization problem with application in portfolio selection}}
\author{
   Anurag Jayswal,
   Sunidhi\textsuperscript{*}\\[5pt]
  \textit{Department of Mathematics and Computing, Indian Institute of Technology}\\
  \textit{(Indian School of Mines), Dhanbad-826004, Jharkhand, India}
  }
\date{}

\begin{document}
\maketitle
\renewcommand{\thefootnote}{*}
\footnotetext{
Corresponding author\\
E-mail addresses: anurag@iitism.ac.in (A. Jayswal), sunidhis2031@gmail.com (Sunidhi)
}

\section*{Abstract}
In this article, we consider a class of unconstrained multiobjective interval-valued optimization problems in which each objective function is a sum of a generalized Hukuhara($gH$)-differentiable interval valued function and a proper, convex and $gH$-lower semicontinuous interval valued function, which is possibly non-$gH$-differentiable. Within this framework, we introduce the notion of a proximal gradient direction and characterize the Pareto critical points of the considered problem through a strongly convex auxiliary subproblem. Based on these characterizations, we develop a proximal gradient algorithm with an Armijo-type line search that ensures sufficient descent with respect to the interval partial order.
We then establish the existence of an uniform lower bound on the admissible step length, the global convergence of the generated sequence to a Pareto critical point followed by a global rate of convergence in the general nonconvex setting. The theoretical findings of the proposed method are validated through some test problems and further illustrate its applicability on portfolio selection problems involving interval-valued data.

\medskip
\noindent \textbf{Keywords:}
Multiobjective interval valued optimization; Pareto critical point; Proximal gradient direction; Proximal gradient method.

\medskip
\noindent \textbf{MSC classification:} 
65K05, 90C29, 90C70

\section{Introduction}

Real world decision making problems arising in engineering design, economics, and finance rarely admit a single criterion of optimality instead, they typically require the simultaneous optimization of several, often conflicting, objective functions. This gives rise to the framework of multiobjective optimization, wherein the classical notion of a unique minimizer is replaced by the notion of Pareto optimality, and the goal shifts to computing a solution that cannot be improved in one objective without deteriorating another. Owing to this rich solution structure, multiobjective optimization problems (MOP) have found extensive applications in portfolio selection, engineering design, machine learning, and resource allocation, among others.

\medskip
\noindent Among the various solution approaches to MOP, descent-type methods that generalize classical scalar-valued unconstrained optimization algorithms have received considerable attention.  Fliege and Svaiter~\cite{fliege2000} proposed the steepest descent method for unconstrained multicriteria optimization, which was later extended in several directions, including the projected gradient method proposed by Drummond and Iusem~\cite{granadrummond2004} for constrained problems, and the Newton method by Fliege $et\ al.$~\cite{fliege2009newton}, which achieves a faster local rate of convergence by incorporating second-order information of the objectives. Further, Bonnel $et\ al.$~\cite{bonnel2005} introduced a proximal point method for nonsmooth vector optimization problems, while Tanabe $et\ al.$~\cite{TanabePGM} developed a proximal gradient method for multiobjective optimization problems having composite objective structure, and later on analyzed its convergence results in Tanabe $et\ al. $ \cite{converegnceTanabe}, extending the scalar proximal gradient framework of Beck and Teboulle~\cite{beck2009} to the Pareto setting.

\medskip
\noindent Parallel to these developments, a substantial body of literature has emerged for handling optimization problems in which the data of the problem are not known precisely but are only available within some range of uncertainty. Interval analysis, initiated by Moore~\cite{moore1966interval}, provides a natural and computationally tractable framework for representing such uncertainty through compact intervals, avoiding both the restrictive worst-case bounds of robust optimization and the distributional assumptions of stochastic optimization. Since the usual interval difference fails to preserve fundamental algebraic properties, Stefanini~\cite{stefanini2008generalization} introduced the generalized Hukuhara ($gH$) difference, which subsequently enabled the development of a full-fledged $gH$-differential calculus for interval-valued functions (IVFs). Building on this calculus and the partial order relations for intervals studied by Wu~\cite{wu2011interval}, optimality conditions of Karush--Kuhn--Tucker type for interval-valued optimization problems (IVOP) were established by Wu~\cite{Wu2007}, while duality and saddle point results were further studied by
Jayswal $et\ al.$\ \cite{Jayswal2011} and Ahmad $et\ al.$\ \cite{Ahmad2016}.

\medskip
\noindent Other than the optimality study, by the success of descent-type algorithms in the scalar multiobjective setting, various iterative methods have also recently been developed for multiobjective interval-valued optimization problems (MIVOP).  Upadhyay $et\ al.$ \cite{UpadhyayNewton} proposed a Newton method to  solve MIVOP and further a generalized Newton approach was developed in Upadhyay $et\ al$ \cite{genralNewtonBB}. Quasi-Newton method was discussed in Upadhyay $et.\ al$ \cite{quasiUpadhyay} to accelerate local convergence near a Pareto critical point. Subsequently, a conjugate direction method was proposed by Upadhyay $et\ al.$ \cite{conjugateUpadhay}. Recently, Mondal and Ghosh~\cite{mondal2025steepest} proposed a steepest descent method for MIVOP to capture the Pareto critical points. However, all of these existing descent methods for MIVOPs are restricted to objective functions that are $gH$-differentiable throughout their domain and consequently, are not applicable when the objective possesses a composite structure involving a non $gH$-differentiable component, such as a sparsity-inducing regularizer or a nonsmooth cost term. To the best of our knowledge, no study has so far been developed for MIVOP with a composite objective structure that involves both the smooth and non-smooth interval valued functions, despite the fact that such composite formulations arise naturally whenever a smooth risk or loss term is combined with a nonsmooth regularization or cost term under interval uncertainty in the problem data. 

\medskip
\noindent Motivated by the above discussion, and building upon the vector proximal gradient framework of Tanabe $et\ al.$~\cite{TanabePGM} and extending the interval-valued steepest descent framework of Mondal and Ghosh~\cite{mondal2025steepest}, we develop a proximal gradient method for a class of multiobjective interval-valued optimization problems with composite objective functions (MIVOP$_C$). The proposed framework extends the existing interval-valued optimization methodology to composite settings and provides a systematic approach to computing Pareto critical points.

\medskip
\noindent The main contributions of this paper are listed below:
\begin{itemize}
\item We establish the relationship between the solution concepts of the MIVOP$_C$. 
    \item We introduce the notion of a proximal gradient direction and characterize the Pareto critical points of MIVOP$_C$ through a strongly convex auxiliary subproblem.
    \item We propose an implementable algorithm for the proximal gradient method which is equipped with an Armijo-type line search to determine the step length. In addition, we provide the existence of a uniform lower bound on the step lengths.
    \item We derive the global convergence of the sequence generated by the proposed method to a Pareto critical point and further, discuss its rate of convergence in the nonconvex setting.
    \item Finally, we illustrate the practical relevance of the proposed framework through some test problems and extend it applicability to solve portfolio selection problem with interval-valued financial parameters.
\end{itemize}

\noindent The remainder of this article is organized as follows. In Section~\ref{section 2}, the basic concepts and the notions of Pareto optimality and Pareto criticality of the multiobjective interval valued optimization problem (MIVOP$_C)$ are  given. Section~\ref{section 3} develops the proximal gradient method for the MIVOP$_C$. Section~\ref{section 4} establishes the convergence analysis of the proposed method. In Section~\ref{Section 5}, we perform numerical experiments of the proposed method on some test problems and further, its applicability to a portfolio selection problem is discussed in Section~\ref{section 6}. Finally, Section~\ref{section 7} concludes the paper and discusses the possible future directions.

\section{Preliminaries}\label{section 2}

In this section, we provide some standard notations used throughout the article and recall some basic concepts related to proximal and interval analysis that form the fundamental tool of the proposed methodology. Finally, we develop the necessary solutions concepts for the considered multiobjective interval valued optimization problem.

\subsection*{Notations and basic concepts}
We use the following notations in this article:  
$\mathbb{R}$, $\R_+$ $\R_{++}$,
$\mathbb{R}^n$, $\R_{\infty}$ = $\mathbb{R} \cup \{ +\infty\}$, $\mathbb{N}$, $\|\cdot\|$, 
$\IR$,
$\IR_{\infty} = \IR \cup \{[ +\infty,+\infty]\}$, $\IR_+$ and $\IR^n$ to represent the real line, the positive real line, the non-zero positive real line, $n$-dimensional Euclidean space, the extended real line, the set of natural numbers, the Euclidean norm, the set of compact intervals defined on the real line, the family of compact intervals defined on the extended real line, the family of compact intervals on non-negative real line and the Cartesian product space of compact intervals $\IR\times \IR\times \dots\times\IR(n$ times$)$, respectively.

\medskip
\noindent The following basic concepts and definitions are from Rockafellar~\cite{tyrrell1970convex}.

\medskip
\noindent Let $\zeta:\R^n\to \R_{\infty}$ be a real valued function that takes values from the extended space, then its effective domain, denoted by $dom\ \zeta$, is given as a set $\{x\in \R^n\ |\ \zeta(x)<+\infty\}$. The function $\zeta$ is said to be proper if $\text{dom}\ \zeta\neq \emptyset$.

\begin{definition}(Convex function ~\cite{tyrrell1970convex}) 
The function $\zeta$ is said to be convex, if for all $x,y\in \text{dom}\ \zeta$ and $\lambda\in [0,1]$,
\[
\zeta(\lambda x+(1-\lambda)y)\le \lambda \zeta(x)+(1-\lambda)\zeta(y).
\]
The function $\zeta$ is said to be strictly convex if the inequality $`\le$' is replaced by $`<$'.
\end{definition}

\begin{definition}(Strongly convex ~\cite{tyrrell1970convex})
The function $\zeta $ is said to be strongly convex with respect to modulus $\nu>0$, if for all $x,y\in \text{dom}\ \zeta$ and $\lambda\in [0,1]$,
\[
\zeta(\lambda x+(1-\lambda)y)\le \lambda \zeta(x)+(1-\lambda)\zeta(y)-\frac{\nu}{2}\lambda(1-\lambda)\|x-y\|^2.
\]
\end{definition}

\begin{definition}\label{lsc}
(Lower semicontinuity ~\cite{tyrrell1970convex})
The function $\zeta$ is said to be lower semicontinuous over $\R^n$, if for every sequence $\{x^k\}\subset \R^n$ with $x^k\to x$, we have 
\[
\zeta (x)\le \lim_{k\to \infty}\inf \zeta(x^k),\ \forall x\in \R^n.
\]
\end{definition}
\noindent Throughout this article, unless stated otherwise, every extended real valued function considered in this paper belongs to the set $\Gamma_0(\R^n)$, defined as
\[
\Gamma_0(\R^n)=\{\zeta:\R^n\to \R_{\infty} \ | \ \zeta \text{ is convex, proper and lower semicontinuous}\}.
\]

\begin{definition}(Directional derivative \cite{bertsekasConves})
For any  $x,d\in\R^n$, the directional derivative of the function $\zeta$ at $x$ along $d$ is denoted by $\zeta'(x;d)$ and is given by
\[
\zeta'(x;d):=\lim_{\lambda\to 0}\frac{\zeta(x+\lambda d)-\zeta(x)}{\lambda}.
\]
If $\zeta$ is differentiable then, the above reduces to $\zeta'(x;d):=\nabla \zeta(x)^\top d$, where $\nabla \zeta(x)$ is the gradient of $\zeta $ at point $x$.
\end{definition}\label{def: dd scalar}
\begin{lemma}\cite{bertsekasConves}
    If $\zeta$ is a convex function, then we have the following non-decreasing property,
\[
\zeta(x+d)-\zeta(x)\ge \frac{\zeta(x+\lambda d)-\zeta(x)}{\lambda},\ \text{for all }\lambda\in(0,1).
\]
\end{lemma}

\begin{definition}\label{subgrad}
(Subgradient ~\cite{tyrrell1970convex})
    Let $\zeta:\R^n\to \R$ be a convex function. Then $\xi$ is called the  subgradient of $\zeta$ at the point $x\in \R^n$, if
    \[
    \zeta(y)-\zeta(x)\ge \langle \xi,y-x\rangle,
    \]
    for all $y\in \R^n$. The collection of all subgradients of the function $\zeta$ at the point $x\in \R^n$ is called the subdifferential of $\zeta$ at $x$, denoted by $\partial \zeta$, and is given as
    \[
    \partial \zeta :=\{\xi\in\R^n\ |\ \zeta(y)-\zeta(x)\ge \langle \xi,y-x\rangle,\ \forall y\in \R^n\}.
    \]
\end{definition}

\subsection*{Interval analysis}
Let $\mathcal U\in \IR$, then $\mathcal U$ can be represented as $\mathcal U=[\underline{\mathcal U},\overline{\mathcal U}]$, where $\underline{\mathcal U}\le \overline{\mathcal U}$ and $\underline{\mathcal U}$,$\overline{\mathcal U}\in \R$.
According to Moore \cite{moore1966interval}, for a pair of intervals
$\mathcal U=[\underline{\mathcal U},\overline{\mathcal U}]$ and $\mathcal V=[\underline{\mathcal V},\overline{\mathcal V}]$ in $\IR$, and a real number $c$, the interval arithmetic are denoted as $\mathcal U\oplus \mathcal V$, $c\odot \mathcal V$, and $\mathcal U\odot \mathcal V$, respectively, and are defined as:
\begin{itemize}
    \item $\mathcal U\oplus \mathcal V=[\underline{\mathcal U}+\underline{\mathcal V},\ \overline{\mathcal U}+\overline{\mathcal V}].$
    \item $ c \odot \mathcal U =
    \begin{cases}
[c \underline{\mathcal U},\ c \overline{\mathcal U}], & c \ge 0, \\
[c \overline{\mathcal U},\ c \underline{\mathcal U}], & c < 0 .
\end{cases}$
   \item $ \mathcal U\odot \mathcal V$= [min$\{\underline{\mathcal U}.\underline{\mathcal V}, \underline{\mathcal U}.\overline{\mathcal V}, \overline{\mathcal U}.\underline{\mathcal V}, \overline{\mathcal U}.\overline{\mathcal V}\}$, max$\{\underline{\mathcal U}.\underline{\mathcal V}, \underline{\mathcal U}.\overline{\mathcal V}, \overline{\mathcal U}.\underline{\mathcal V}, \overline{\mathcal U}.\overline{\mathcal V}\}].$
\end{itemize}
In addition, if $\mathcal U\in \IR$ such that $\underline{\mathcal U}=\overline{\mathcal U}$, then we call $\mathcal U$ as a degenerate interval. In this case, we denote
\[
\mathbf{\mathcal U}=[a,a],\quad a\in \R,
\]
where $a=\underline{\mathcal U}=\overline{\mathcal U}$. In particular, the zero interval is denoted by $\textbf{0}=[0,0]$.\\
\\
\noindent The usual difference of two compact intervals is not suitable for interval calculus, since it does not preserve basic algebraic properties (see Moore \cite{moore1966interval}). To overcome this challenge, Stefanini \cite{stefanini2008generalization} introduced the concept of generalized Hukuhara ($gH$) difference.

\begin{definition}\label{def: gh difference} 
($gH$-difference~\cite{stefanini2008generalization})
Let $\mathcal U,\mathcal V$ and $\mathcal W$ be three arbitrary elements of $\IR$ such that $\mathcal U=\mathcal V\oplus \mathcal W$ or $\mathcal V=\mathcal U\ominus \mathcal W$, then $\mathcal W$ is said to be the \textit{$gH$-difference} between $\mathcal U$ and $\mathcal V$, and it is denoted as $\mathcal W=\mathcal U\ominus_{gH} \mathcal V$. In other words, the $gH$-difference between $\mathcal U$ and $\mathcal V$ is defined as:
\[
   \mathcal U\ominus_{gH}\mathcal V=[\min\{\underline{\mathcal U}-\underline{\mathcal V},\ \overline{\mathcal U}-\overline{\mathcal V}\},\ \max\{\underline{\mathcal U}-\underline{\mathcal V},\ \overline{\mathcal U}-\overline{\mathcal V}\}].
\]
\end{definition}

\noindent In interval analysis, order relations are used to compare the interval values and thus, several partial orderings of the intervals have been discussed in the literature (see \cite{moore1966interval, mondal2025steepest, wu2011interval}, and the references cited therein). We consider the partial ordering between the intervals $\mathcal U,\mathcal V \in \IR$ from Wu \cite{wu2011interval} as :
\[
\mathcal U \preceq_{LU} \mathcal V \quad \Leftrightarrow \quad \underline{\mathcal U}\le \underline{\mathcal V}
\ \text{and}\ 
\overline{\mathcal U}\le \overline{\mathcal V},
\]
where "$\preceq_{LU}$" refers to the partial ordering in $\IR$. Thus, $\mathcal U$ is inferior to $\mathcal V$ if $\mathcal U$ lies entirely inside $\mathcal V$. Also, $\mathcal U\prec_{LU}\mathcal V$ if and only if $\mathcal U\preceq_{LU}\mathcal V$ and $\mathcal U\neq \mathcal V$ which means $\mathcal U\prec_{LU}\mathcal V$ if and only if :
\[
\begin{cases}
    \underline{\mathcal U}< \underline{\mathcal V}\\
    \overline{\mathcal U}=\overline{\mathcal V}
\end{cases}
,\quad \text{or} \quad
\begin{cases}
    \underline{\mathcal U}= \underline{\mathcal V}\\
   \overline{\mathcal U}<\overline{\mathcal V}
\end{cases}
,\quad \text{or} \quad\
\begin{cases}
    \underline{\mathcal U}< \underline{\mathcal V}\\
    \overline{\mathcal U}<\overline{\mathcal V}
\end{cases}
.
\]
We also write $\mathcal U\precprec_{LU} \mathcal V$ if and only if $ \underline{\mathcal U}< \underline{\mathcal V}$ and $ \overline{\mathcal U}<\overline{\mathcal V}$. 
For simplicity, we use the notations `$\preceq$' and `$\prec$' instead of `$\preceq_{LU}$' and `$\prec_{LU}$', respectively, in this article.

 \begin{remark}
      For any two arbitrary elements $\hat {\mathcal U}=(\mathcal{U}_1,\mathcal{U}_2,\dots,\mathcal U_n)^\top$ and $\hat {\mathcal V}=(\mathcal V_1,\mathcal V_2,\dots,\mathcal V_n)^\top$ in $\IR^n$, we have
      \begin{align*}
      \hat {\mathcal U}\preceq \hat {\mathcal V}\ &\Leftrightarrow\ \mathcal U_i\preceq \mathcal V_i,\\ 
      \hat {\mathcal U}\prec \hat {\mathcal V}\ &\Leftrightarrow\ \mathcal U_i\prec \mathcal V_i,\ \forall i=1,2,\dots,n.
      \end{align*}
 \end{remark}

\begin{definition}\label{IVF norm}
(Norm on $\IR$ and $\IR^n$~\cite{moore1966interval})
Let $\mathcal U=[\underline{\mathcal U},\overline{\mathcal U}]\in \IR$, then the norm on $\IR$ is given by
\[
   \|\mathcal U\|_{\IR}=max\{|\underline{\mathcal U}|,|\overline{\mathcal U}|\}.
\] 
Further, the norm on $\IR^n$, denoted as $\|\cdot\|_{\IR^n}$, is given by
\[
\|\hat {\mathcal U}\|_{\IR^n}=\sum_{i=1}^n\|\hat {\mathcal U}_i\|_{\IR},
\]
where $\hat {\mathcal U}:=(\hat {\mathcal U}_1,\hat {\mathcal U}_2,\dots,\hat {\mathcal U}_n)\in \IR^n$ is an interval valued vector.
\end{definition}

\noindent The following are few properties of intervals:

\begin{lemma}\label{lemma ppt}
   \cite{Ghosh2021} Let $\mathcal U,\mathcal V,\mathcal Y,\mathcal Z\in \IR$ be arbitrary elements. Then, the following results hold:
    \begin{enumerate}[label=(\roman*)]
    \normalfont
        \item \textit{If} $\mathcal U\preceq \mathcal V$ \textit{and} $\mathcal V\preceq \mathcal Y$, \textit{then} $\mathcal U\preceq \mathcal Y$,
        \item \textit{If} $\mathcal U\oplus \mathcal V\preceq \mathcal Y\oplus \mathcal Z$, \textit{then} $\mathcal U\ominus_{gH}\mathcal Y\preceq \mathcal Z\ominus_{gH}\mathcal V $.
    \end{enumerate}
\end{lemma}

\begin{remark}\label{remark 1}
\begin{enumerate}[label=(\roman*)]
   \item If $\mathcal V=\mathbf{0}$, then $\mathcal U\preceq \mathcal Y\oplus \mathcal Z\ \Rightarrow \mathcal U\ominus_{gH}\mathcal Y\preceq \mathcal Z$.
   \item If $\mathcal U\preceq \mathcal V$, then $\mathcal U\oplus \mathcal Y\preceq  \mathcal V\oplus \mathcal Z$, then from Lemma~\ref{lemma ppt}, $\mathcal U\ominus_{gH} \mathcal Y\preceq \mathcal V\ominus_{gH}\mathcal Y$ or $\mathcal Y\ominus_{gH} \mathcal V\preceq \mathcal Y\ominus_{gH}\mathcal U$.
 \end{enumerate}  
\end{remark}

\begin{proposition}\label{inclusion lemma}
  \cite{Hai2026}  If $\mathcal U,\mathcal V,\mathcal Y,\mathcal Z\in \IR$ are arbitrary elements, then the following result holds:
    \[
    (\mathcal U\oplus \mathcal Y)\ominus_{gH}(\mathcal V\oplus \mathcal Z)=(\mathcal U\ominus_{gH}\mathcal V)\oplus (\mathcal Y\ominus_{gH}\mathcal Z),
    \]
    if either of the conditions $\overline{\mathcal U}-\underline{\mathcal U}\ge \overline{\mathcal V}-\underline{\mathcal V}$ and $\overline{\mathcal Y}-\underline{\mathcal Y}\ge \overline{\mathcal Z}-\underline{\mathcal Z}$ or $\overline{\mathcal U}-\underline{\mathcal U}\le \overline{\mathcal V}-\underline{\mathcal V}$ and $\overline{\mathcal Y}-\underline{\mathcal Y}\le \overline{\mathcal Z}-\underline{\mathcal Z}$ satisfies.
\end{proposition}

\begin{remark}  The width conditions required for equality in the Proposition~\ref{inclusion lemma} are not restrictive in many practical settings, as they are naturally satisfied by interval valued mappings with monotonic or constant widths. Consequently, the resulting equality can be used conveniently in the subsequent analysis. (See, for instance \cite{Hai2026})
\end{remark}

\subsection*{Interval valued functions}
A function $F : \R^n \to \IR$ is called an interval valued function (IVF), represented as $F(x):= [\underline{F}(x), \overline{F}(x)]$, where the endpoint functions $\underline{F}, \overline{F} : \R^n \to \R$ satisfy
$\underline{F}(x) \le \overline{F}(x), \ \forall x \in \R^n$. If $F:\R^n\to \IR_{\infty}$, then it is called an extended IVF. In this setting, the endpoint functions of $F$ take values in the extended real line $\R_{\infty}$. 

\medskip
\noindent Throughout the article, all the interval valued functions considered are assumed to satisfy the conditions of Proposition~\ref{inclusion lemma}.

\begin{definition}($gH$-continuity~\cite{mondal2025steepest}) An IVF $F:\R^n\to \IR$ is said to be $gH$-continuous at a point $\hat x\in \R^n$ if 
\[
\lim_{\|\lambda\|\to 0}F(\hat x+\lambda)\ominus_{gH}F(\hat x)=\mathbf{0}.
\]
Further, the IVF $F$ is said to be $gH$-continuous on $\R^n$, if the above condition holds for every $x\in \R^n$.
\end{definition}

\begin{definition}($gH$-Lipschitz continuous~\cite{mondal2025steepest}) An IVF $F:\R^n\to \IR$ is said to be $gH$-Lipshcitz continuous if there exists a constant $L>0$ such that for all $x,y\in \R^n$,
\[
\|F(x)\ominus_{gH}F(y)\|_{\IR}\le L\|x-y\|.
\]
\end{definition}

\begin{lemma}\label{lemma gh lip}
\cite{mondal2025steepest} The IVF $F:\R^n\to \IR$ is $gH$-continuous ($gH$-Lipschitz continuous) on $\R^n$ if and only if the endpoint functions $\underline{F}$ and $\overline{F}$ are continuous (Lipschitz continuous) on $\R^n$.
\end{lemma}

\begin{definition}\label{eff dom IVF}
(Effective domain of IVF \cite{kumar2023ekeland})
The set of all the points at which an extended IVF $F:\R^n\to \IR_{\infty}$, takes a finite interval value is called its effective domain and is expressed as
\[
dom\ F:=\{x\in \R^n\ |\ F(x)\prec\mathbf{+\infty}\}.
\]
\end{definition}

\noindent Note that, by the order relations defined above, it is obvious that $$dom\ F:= \text{dom}\ \underline{F}\cap \text{dom}\ \overline{F}\ ,$$
where $\underline{F},\ \overline{F}$ are endpoint functions of the IVF $F$.

\begin{definition}(Proper IVF \cite{kumar2023ekeland})  An extended IVF $F:\R^n\to \IR_{\infty}$ is said to be proper if its effective domain is non-empty, $i.e.,$ $dom\ F\neq \emptyset$.
\end{definition}

\begin{definition}(Convex IVF \cite{kumar2023ekeland}) An extended IVF $F:\R^n\to \IR_{\infty}$ is said to be convex if for all $x,y\in dom\ F$ and $\hat \lambda\in[0,1]$, one has 
\[
F(\hat \lambda x+(1-\hat \lambda)y)\preceq \hat \lambda F(x)+(1-\hat \lambda )F(y).
\]
\end{definition}

\begin{lemma}\label{convex IVF lemma} 
\cite{Wu2007}An IVF $F:\R^n\to \IR$ is said to be convex if and only if the endpoint functions $\underline{F}$ and $\overline{F}$ are convex.
\end{lemma}

\begin{lemma}\label{gi gh-continuous}
    Let $F:\R^n\to \IR_{\infty}$ be a proper convex IVF, and let $\mathcal W$ be a compact subset of $int(dom\ F)$, then $F$ is $gH$-Lipschitz continuous on the set $\mathcal W$, provided $int(dom\ F):= int(dom\ \underline{F})\cap int(dom\ \overline{F}) \neq \emptyset$.
\end{lemma}
\begin{proof}
Since the IVF $F$ is proper convex, the endpoint functions are proper convex, therefore from Theorem 10.4 of Rockafellar \cite{tyrrell1970convex}, the end point functions of $F$ are Lipschitz continuous on $\mathcal W$, and hence from Lemma~\ref{lemma gh lip}, the IVF $F$ is $gH$-Lipschitz continuous on the set $\mathcal W$.       
\end{proof}

\begin{definition}($gH$-derivative~\cite{Kumar2022weak}) Let $\mathcal V\subseteq \R$ be a nonempty open set. The $gH$-derivative of an IVF $F:\mathcal V\to \IR$ at a point $\hat x\in \mathcal V$, denoted by $F'(\hat x)$, is defined as
\[
F'(\hat x)=\lim_{d\to 0}\frac{1}{d}\odot (F(\hat x+d)\ominus_{gH}F(\hat x)),
\] 
provided the limit exists.
\end{definition}

\begin{definition}\label{def: gh dd}
(gH-directional derivative~\cite{Kumar2022weak})
Let $\mathcal W\subseteq \R^n$ be a nonempty set. The $gH$-directional derivative of an IVF $F:\mathcal W\to \IR$ at a point $\hat x\in \mathcal W$ in the direction $d\in \R^n$, denoted by $F'_{gH}(x;d)$, is defined as
\[
F'_{gH}(\hat x;d)=\lim_{\lambda\to 0^+}\frac{1}{\lambda}\odot F(\hat x+\lambda d)\ominus_{gH}F(\hat x),
\]
provided the limit exists.
\end{definition}

\begin{lemma}\label{DD Lemma}
\cite{Kumar2022weak}
    Let $F:\R^n\to \IR_{\infty}$ be a proper convex IVF. Then, it follows that
    \[
    \lim_{\lambda\to 0^+}\frac{1}{\lambda}\odot F(x+\lambda d)\ominus_{gH}F(x)\preceq F(x+d)\ominus_{gH}F(x),
    \]
     where $d=y-x$, for all $x,y\in dom \ F$.
\end{lemma}

\begin{definition}($gH$-gradient \cite{mondal2025steepest}) Let $\mathcal W\subseteq \R^n$ be a nonempty open set. The $gH$-gradient of an IVF $F:\mathcal W \to \IR$ at a point $\hat x=(\hat x_1,\hat x_2,\dots,\hat x_n)\in \R^n$, is defined as
\[
\nabla_{gH}F(\hat x)=(F'_1(\hat x_1),F'_2(\hat x_2),\dots,F'_n(\hat x_n))^\top,
\]
where $F'_j(\hat x_j)$ is the $gH$-derivative of the IVF $F_j(x_j)= F(\hat x_1,\hat x_2,\dots,\hat x_{j-1},x_j,\hat x_{j+1},\dots,\hat x_n)$ at $x_j=\hat x_j$, for every $j=\overline{1,n}$.  
\end{definition}

\begin{definition}($gH$-Lipshictz gradient \cite{mondal2025steepest}) An IVF $F:\mathcal W\to \IR$ is said to possess $gH$-Lipschitz gradient on a nonempty open set $\mathcal W\subseteq \R^n$, if there exists a constant $L>0$ satisfying
\[
\|\nabla_{gH}F(x)\ominus_{gH}\nabla_{gH}F(y)\|_{\IR^n}\le L\|x-y\|,\ \forall x,y\in \mathcal W.
\]    
\end{definition}

\begin{definition}\label{gh-diff def}
($gH$-differentiability~\cite{Kumar2022weak}) Let $\mathcal W\subseteq \R^n$ be a nonempty open set. An IVF $F:\mathcal W\to \IR$ is said to be $gH$-differentiable at a point $\hat x\in \mathcal W$, if there exists two IVFs $\mathcal P(F(\hat x);d)$ and $\mathcal Q(\hat x;d):\mathcal \R^n\to \IR$, and a constant $\delta>0$ such that 
\[
F(\hat x+d)\ominus_{gH} F(\hat x)=\mathcal Q(\hat x;d)\oplus \|d\|\odot \mathcal P(F(\hat x);d),
\]
for every $d$ satisfying $\|d\|<\delta$, where $\mathcal P(F(\hat x);d)\to 0$ as $\|d\|\to 0$ and $\mathcal Q(\hat x;d)$ satisfies the following properties,
\begin{enumerate}[label=(\roman*)]
    \item $\mathcal Q(\hat x;d_1+d_2)=\mathcal Q(\hat x;d_1)\oplus \mathcal Q(\hat x;d_2),\ \forall d_1,d_2\in \mathcal W$.
    \item $\mathcal Q(\hat x;p\hat d)=p\odot \mathcal Q(\hat x;\hat d),\ \forall p\in \R$ and $\hat d\in \mathcal W$. 
\end{enumerate}
\end{definition}

\begin{lemma}\label{gh diff lemma}
\cite{Kumar2022weak}
     Let an IVF $F:\mathcal W\to \IR$ be $gH$-differentiable at a point $\hat x\in \mathcal W$, where $\mathcal W\subseteq\R^n$ is a nonempty open set. Then there exist a non zero $\lambda$ such that 
     \[
     \lim_{\lambda\to 0^+}\frac{1}{\lambda}\odot F(\hat x+\lambda d)\ominus_{gH}F(\hat x) =\mathcal Q(\hat x;d),
     \]
    for all $d\in \R^n$ with $\|\lambda d\|<\delta$, where $\mathcal Q(\hat x;d)$ is an IVF as mentioned in Definition~\ref{gh diff lemma}. Further, upon the existence of the $gH$-gradient of the IVF $F$ at the point $\hat x\in \mathcal W$, the IVF $\mathcal Q(\hat x;d)$ is termed as
    \[
    \mathcal Q(\hat x;d)= \nabla_{gH}F(\hat x)^\top\odot d,
    \]
    for all $d=(d_1,d_2,\dots,d_n)^\top\in \R^n$.
\end{lemma}

\begin{definition}(Infimum and supremum of an extended IVF~\cite{kumar2023ekeland}) For a nonempty set $\mathcal W\subseteq\R^n$, consider an extended IVF $F:\mathcal W\to \IR_{\infty}$. The infimum of $F$ over $\mathcal W$ is defined as the infimum of its range set, namely
\[
\inf_{x \in \mathcal W} F(x)
:=
\inf\{F(x):x\in\mathcal{W}\}.
\]
In the analogous manner, the supremum of $F$ over $\mathcal W$ is given by
\[
\sup_{x \in \mathcal W} F(x)
:=
\sup\{F(x):x\in\mathcal{W}\}.
\]
\end{definition}

\begin{definition}\label{def:gHLsc}
($gH$-Lower Semicontinuity of an extended IVF \cite{kumar2023ekeland}) Let $F:\R^n \to \IR_{\infty}$ be an extended IVF and let $x^*\in \R^n$. For $\delta>0$, we define the open ball $B_{\delta}(x^*):=\{x\in \R^n: \|x-x^*\|<\delta\},$ and the limit inferior of $F$ at $x^*$ as
\[
\lim_{x\to x ^*}\inf F(x)=\sup_{\delta >0}\inf \{F(x): x\in B_\delta (x^*)\}.
\]
Then, the IVF $F$ is a generalized Hukuhara lower semicontinuous ($gH$-$l.sc$) at $x^*\in \R^n$, if
\[
F(x^*)\preceq_{LU}\lim_{x\to x ^*}\inf F(x),\quad \forall x\in B_{\delta}(x^*).
\]
Further, the IVF $F$ is said to be $gH$-$l.sc$ on $\R^n$, if the above inequality holds for every $x^*\in \R^n$.
\end{definition}
\begin{remark}\label{Remark 2.5} In view of the Definition \ref{def:gHLsc}, it follows that the extended IVF $F$ is $gH$-lower semicontinuous at the point $x^*\in \R^n$ if and only if the end point functions are lower semicontinuous at the point $x^*\in \R^n$ in the sense of Definition \ref{lsc}. Furthermore, if an extended IVF $F$ is $gH$-continuous on $\R^n$, then $F$ is $gH$-lower semicontinuous on $\R^n$ (see, for instance, Kumar and Ghosh \cite{kumar2023ekeland}).
\end{remark}

\subsection*{Multiobjective interval valued composite optimization problem}

We consider the following unconstrained multiobjective interval valued composite optimization problem 
\begin{align*}
    \textbf{MIVOP}_C\qquad \min_{x\in\R^n} F(x):=(F_1(x),\dots, F_m(x))^\top,
\end{align*}
where each component $F_i:\R^n\to \IR_{\infty}$ admits an interval valued composite structure, defined as 
\[
F_i(x)=f_i(x)\oplus g_i(x),\quad \forall i=\overline{1,m}.
\]
Here, the IVF $f_i:\R^n\to \IR$ is $gH$-differentiable with $gH$-Lipschitz continuous gradients and the IVF $g_i:\R^n\to \IR_{\infty}$ is convex, proper and $gH$-lower semicontinuous, not necessarily $gH$-differentiable, whose endpoint functions $\underline{g_i},\overline{g_i}\in \Gamma_0(\R^n)$, for each $i=\overline{1,m}$.

\medskip
\noindent The following are the definitions of the solution concepts of the MIVOP$_C$. 

\begin{definition}[Weak Pareto optimality] 
A point $\hat x\in \R^n$ is said to satisfy weak Pareto optimality for MIVOP$_C$ if no feasible point $x\in \R^n$ exists such that, for all $i=\overline{1,m}$,
\[
F_i(x)\prec F_i(\hat x).
\]
Then $\hat x$ is called a weakly Pareto point of MIVOP$_C$.
\end{definition}

\begin{definition}[Pareto optimality] 
A point $\hat x\in \R^n$ is regarded as a Pareto optimal point of MIVOP$_C$ whenever no other point $x\in \R^n$ exists such that, for all $i=\overline{1,m}$,  
\[
F_i(x)\preceq F_i(\hat x).
\]
\end{definition}

\begin{definition}[Pareto criticality]\label{PCP}
A point $\hat x\in \R^n$ is termed as Pareto critical if there exist no such $d\in\R^n$ that satisfy
\begin{align}
    F'_{gH,i}(x;d)\prec \mathbf{0},\ \forall i=\overline{1,m}, \label{pcp condition}
\end{align}
where $F'_{gH,i}(x;d)$ is $gH$-directional derivative of each of the components of the objective function $F(x)=(F_1(x),\dots,F_m(x))^\top$ of the MIVOP$_C$.

\medskip
\noindent Since the IVF $f_i$ is $gH$-differentiable, using the Lemma~\ref{gh diff lemma}, the above condition~\eqref{pcp condition} can be equivalently written as
\begin{align}
    \nabla_{gH}f_i(x)^\top\odot d \oplus g'_{gH,i}(x;d)\prec \mathbf{0},\ \forall i=\overline{1,m}.
\end{align}
\end{definition}
\noindent The following result characterizes the relationship between the Pareto criticality and weak Pareto optimality of the MIVOP$_C$.

\begin{proposition}\label{ppsn 1}
\begin{enumerate}[label=(\roman*)]
\normalfont
   \item \textit{Every weakly Pareto optimal point of the MIVOP$_C$ is necessarily a Pareto critical point.}
   \item \textit{If $f_1, f_2,\dots, f_m$ are convex IVFs, then every Pareto critical point of the MIVOP$_C$ is a weakly Pareto optimal point}. 
   \end{enumerate}
\end{proposition}
\begin{proof}
\begin{enumerate}[label=(\roman*)]
\item Let $\hat x\in \R^n$ be a weakly Pareto optimal point of the MIVOP$_C$. Proceeding on contrary, we assume that $\hat x$ is not a  Pareto critical point, therefore by Definition~\ref{PCP}, there exists $d\in \R^n$ satisfying 
\begin{align*}
&F^{'}_{gH,i}(\hat x;d)\prec \mathbf{0}, \ \forall i=\overline{1,m}, \\
\text{or, }\lim_{\lambda\to 0^+}&\frac{1}{\lambda}\odot \big(F_i(\hat x+\lambda  d)\ominus_{gH} F_i(\hat x)\big )\prec \mathbf{0}, \ \forall i=\overline{1,m}.
\end{align*}
Consequently, there exists some $\delta>0$ such that for every $\lambda\in(0,\delta)$,
\begin{align*}
  F_i(\hat x+\lambda d)\ominus_{gH}F_i(\hat x)&\prec \mathbf{0},\\
  \text{or, } F_i(\hat x+\lambda d)\prec F_i(\hat x),
\end{align*}
for all $i=\overline{1,m}$, which contradicts the assumption that $\hat x$ is weakly Pareto optimal point of the MIVOP$_C$, which concludes the proof.

\item Suppose that $\hat x\in \R^n$ is a Pareto critical point but not weakly Pareto optimal point of the MIVOP$_C$. Then there exists some $\overline x(\neq \hat x)\in \R^n$ such that
\begin{align*}
&F_i(\overline x)\prec F_i(\hat x),\ \forall i=\overline{1,m},\\
\text{or, }&F_i(\overline x)\ominus_{gH}F_i(\hat x)\prec \mathbf{0}.
\end{align*}
Since the IVF $F_i$ is convex, Lemma~\ref{DD Lemma} yields
\begin{align*}
\lim_{\lambda\to 0^+}\frac{1}{\lambda}\odot \big (F_i(\hat x+\lambda(\overline{x}-\hat x))\ominus_{gH}F_i(\hat x)\big)&\preceq F_i(\hat x+(\overline{x}-\hat x))\ominus_{gH}F_i(\hat x),\\
&=F_i(\overline x)\ominus_{gH}F_i(\hat{x}),\\
&\prec \mathbf{0},\quad \forall i=\overline{1,m},
\end{align*}
and therefore, there exists $d=(\overline{x}-\hat x)\in \R^n$ such that $F'_{gH,i}(x;d)\prec \mathbf 0$, for all $i=\overline{1,m}$, which contradicts the assumed Pareto criticality of $\hat x$. This completes the proof.
 \end{enumerate}
\end{proof}

\section{Proximal gradient method}\label{section 3}

In this section, we propose a proximal gradient method for obtaining the Pareto critical points of the MIVOP$_C$. 
We introduce an auxiliary IVF which gives a sufficient condition for the Pareto criticality of MIVOP$_C$. Using it, we construct a strongly convex subproblem whose unique minimizer is the proximal gradient direction, and reformulate it as the tractable constrained problem. We then characterize Pareto critical points through the optimal solution and optimal value of this subproblem. Finally, we develop the computation of a suitable step size together with the summary of the proposed algorithm. 

\medskip
\noindent Since each of the components of the objective IVF $F$ of the  MIVOP$_C$ admits a composite structure, thus for all $d\in \R^n$ and a given point $x\in dom\ F$, we define an IVF $h_i:\R^n\to \IR$ as follows
\begin{align*}
h_i(x;d)&=F_i(x+d)\ominus_{gH}F_i(x),\\
&=(f_i(x+d)\oplus g_i(x+d))\ominus_{gH}(f_i(x)\oplus g_i(x)),\\
&=(f_i(x+d)\ominus_{gH}f_i(x))\oplus(g_i(x+d)\ominus_{gH}g_i(x)),\ (by\ Proposition~\ref{inclusion lemma})\\
&=\nabla_{gH}f_i(x)^\top\odot d\oplus (g_i(x+d)\ominus_{gH}g_i(x)),\\
&=\nabla_{gH}f_i(x)^\top \odot d\ \oplus\Delta g_i(x;d),
\end{align*}
 where $\nabla_{gH}f_i(x)^\top\odot \ d$ is the first order approximation of the $gH$-differentiable IVF $f_i$ and $\Delta g_i(x;d)= g_i(x+d)\ominus_{gH} g_i(x)$, for all $i=\overline{1,m}$.

\begin{proposition}\label{relation hi-F'gh}
    For a given point $x\in dom\ F$ and for each $i=\overline{1,m}$, the following inequality holds
    \[
    F'_{gH,i}(x;d)\preceq h_i(x;d),\ \forall d\in \R^n.
    \]
\end{proposition}
\begin{proof}
    Let $x\in dom\ F$ be given, then by Definition~\ref{def: gh dd}, for each $i=\overline{1,m}$, we have
    \begin{align*}
         F'_{gH,i}(x;d)&=\lim_{\lambda\to 0^+}\frac{1}{\lambda}\odot (F_i(x+\lambda d)\ominus_{gH}F_i(x)),\\
         &=\lim_{\lambda\to 0^+}\frac{1}{\lambda}\odot\Big [ \big (f_i(x+\lambda d)  \oplus g_i(x+\lambda d)\big )\ominus_{gH}\big (f_i(x)\oplus g_i(x)\big )\Big ],\\
         &=\lim_{\lambda\to 0^+}\frac{1}{\lambda}\odot\Big [ \big (f_i(x+\lambda d)\ominus_{gH}f_i(x)\big )\oplus \big(g_i(x+\lambda d)\ominus_{gH}g_i(x)\big) \Big ], \\
         &=\lim_{\lambda\to 0^+}\frac{1}{\lambda}\odot \big (f_i(x+\lambda d)\ominus_{gH}(f_i(x)\big )\oplus\lim_{\lambda\to 0^+}\frac{1}{\lambda}\odot \big (g_i(x+\lambda d)\ominus_{gH}(g_i(x)\big ),\\
         &\preceq \nabla_{gH}f_i(x)^\top \odot d\ \oplus \big (g_i(x+d)\ominus_{gH} g_i(x)\big ),\ (using\ Lemmas\  \ref{DD Lemma}\ and\ \ref{gh diff lemma})\\
         &=h_i(x;d),\ \forall d\in\R^n.
    \end{align*}
    Hence the proof.
\end{proof} 

\begin{corollary}\label{hi<0 impls not PCP}
Consider an arbitrary point $x\in dom\ F$. If a direction $d\in \R^n$ exists such that $h_i(x;d)\prec\mathbf{0}$, then $x$ is not a Pareto critical point of the MIVOP$_C$.
\end{corollary}

\begin{proof}
    The proof follows directly from the Proposition~\ref{relation hi-F'gh} and hence being omitted.
\end{proof}

\begin{corollary}\label{PCP implies no hi<0}
    If $x\in dom\ F$ is a Pareto critical point of the MIVOP$_C$, then there exists no $d\in \R^n$ such that $h_i(x;d)\prec \mathbf{0}$, for all $i=\overline{1,m}$.
\end{corollary}
\begin{proof}
    Suppose that $x\in dom\ F$ is a Pareto critical point and on a contrary, there exists $d\in \R^n$ such that $h_i(x;d)\prec \mathbf{0}$, for all $i=\overline{1,m}$. Then from Corollary~\ref{hi<0 impls not PCP}, $F_{gH,i}(x;d)\prec \mathbf{0}$ which contradicts the fact that $x$ is a Pareto critical point. This completes the proof.
\end{proof}

\begin{proposition}\label{NPC implies hi}
    Suppose $x\in dom\  F$ fails to be a Pareto critical point of the MIVOP$_C$, then there exists $d\in \R^n$ such that $h_i(x; d)\prec \mathbf{0}$, for all $i=\overline{1,m}$.
\end{proposition}
\begin{proof} 
     Suppose $x\in dom\ F$ is not a Pareto critical point, then there exists a vector $\hat d\in \R^n$ such that $F'_{gH,i}(x;\hat d)\prec \mathbf{0}$, for all $i=\overline{1,m}$. Now, for any $s\in \R_{++}$, we have
    \begin{align*}
    \frac{1}{s}\odot h_i(x;s\hat d)
    =\nabla_{gH}f_i(x)^\top\odot \hat d\oplus \frac{1}{s}\odot \big( g_i(x+s\hat d)\ominus_{gH}g_i(x)\big ),
    \end{align*}
    On taking limit $s\to 0^+$, we obtain
\begin{align*}
    \lim_{s\to 0^+}\frac{1}{s}\odot h_i(x;s\hat d)&=
\nabla_{gH}f_i(x)^\top\odot \hat d\oplus g'_{gH,i}(x;\hat d),\\
    &=F'_{gH,i}(x;\hat d)\prec \mathbf{0}.
\end{align*}
Consequently, there exists $\overline{s}>0$, such that for all $s\in (0,\overline{s})$, we have
\[
h_i(x;s\hat d)\prec \mathbf{0}.
\]
On setting $d=s\hat d$, the desired result follows.
\end{proof}
\begin{remark}
It is observable that for a given point $x\in dom\ F$, the IVF $h_i(x;d)$ acts as an auxiliary function for characterizing the Pareto criticality of the point $x$ for the MIVOP$_C$.
\end{remark}


\noindent To work with the IVF $h_i$ computationally, we now write its endpoint functions explicitly as follows
\begin{align*}
\underline{h_i}(x;d)&=\frac{1}{2}\big (  \underline{\nabla_{gH}f_i}(x)+\overline{\nabla_{gH}f_i}(x)\big )^\top d\ -\frac{1}{2}\big( \overline{\nabla_{gH}f_i}(x)- \underline{\nabla_{gH}f_i}(x)\big )^\top|d|+\underline{\Delta g_i}(x;d),\\[4pt]
\overline{h_i}(x;d)&=\frac{1}{2}\big (  \underline{\nabla_{gH}f_i}(x)+\overline{\nabla_{gH}f_i}(x)\big )^\top d\ +\frac{1}{2}\big( \overline{\nabla_{gH}f_i}(x)- \underline{\nabla_{gH}f_i}(x)\big )^\top|d|+\overline{\Delta g_i}(x;d),
\end{align*}
where $|d|=(|d_1|,\dots,|d_n|)^{\top}$, and $\underline{\Delta g_i}$, $\overline{\Delta g_i}$ are the endpoints of the IVF $\Delta g_i(x;d)=g_i(x+d)\ominus_{gH}g_i(x)$.
Note that, the $gH$-gradient part of these endpoint functions follows from the interval arithmetic of the component $\nabla_{gH}f_i(x)^\top\odot d$. (See, for instance Mondal and Ghosh~\cite{mondal2025steepest})

\medskip
\noindent For a given $x\in dom\ F$, let $\psi(x;d):\R^n\to \R$ be a real valued function, defined as
\begin{align}
\psi(x;d):=\max_{i=1,\dots,m}\overline{h_i}(x;d),\ \forall d\in \R^n.
\label{psi def}
\end{align}
Also, since $\underline{h_i}(x;d)\le\overline{h_i}(x;d)$ for each $i=\overline{1,m}$, one has
\begin{align}
    \psi(x;d)<0\ &\Rightarrow h_i(x;d)\prec\mathbf{0},\ \forall d\in \R^n.
    \label{implies}
\end{align}

\begin{lemma}\label{Convex lemma}
    For a given $x\in dom\ F$, the real valued function $\psi(x;d
    ):\R^n\to \R$ is convex and $\psi(x;0)=0$.
\end{lemma}
\begin{proof}
Let $d_\lambda=\lambda d_1+(1-\lambda)d_2$, where $d_1,d_2\in \R^n$ and $\lambda\in [0,1]$, then the expression
\begin{equation}
    \begin{aligned}
        &\frac{1}{2}\big (  \underline{\nabla_{gH}f_i}(x)+\overline{\nabla_{gH}f_i}(x)\big )^\top d_\lambda+\frac{1}{2}\big( \overline{\nabla_{gH}f_i}(x)- \underline{\nabla_{gH}f_i}(x)\big )^\top|d_\lambda|\\[4pt]
        &\quad \quad=\frac{1}{2}\big (  \underline{\nabla_{gH}f_i}(x)+\overline{\nabla_{gH}f_i}(x)\big )^\top (\lambda d_1+(1-\lambda)d_2)\\[4pt]
        &\qquad \hspace{2.5cm}+\frac{1}{2}\big( \overline{\nabla_{gH}f_i}(x)-
         \underline{\nabla_{gH}f_i}(x)\big )^\top|\lambda d_1+(1-\lambda)d_2|,\\[2pt]
        &\quad \quad \le \lambda\Big ( \frac{1}{2}\big (  \underline{\nabla_{gH}f_i}(x)+\overline{\nabla_{gH}f_i}(x)\big )^\top d_1+\frac{1}{2}\big( \overline{\nabla_{gH}f_i}(x)- \underline{\nabla_{gH}f_i}(x)\big )^\top|d_1|\Big )\\[4pt]
        &\quad \quad+(1-\lambda)\Big ( \frac{1}{2}\big (  \underline{\nabla_{gH}f_i}(x) +\overline{\nabla_{gH}f_i}(x)\big )^\top d_2
    +\frac{1}{2}\big( \overline{\nabla_{gH}f_i}(x)- \underline{\nabla_{gH}f_i}(x)\big )^\top|d_2|\Big ).
    \end{aligned}\label{f convex}
\end{equation}
\noindent On the other hand, for a given $x\in dom\ F$, $\overline{\Delta g_i}(x;d)=\overline{g_i(x+d)\ominus_{gH}g_i(x)}$, therefore, one has
\[
\overline{\Delta g_i}(x;d)=\max \{\underline{g_i}(x+d)-\underline{g_i}(x), \overline{g_i}(x+d)-\overline{g_i}(x)\}, \ \forall d\in \R^n.
\]
Since the IVF $g_i$ is convex and the map $d\mapsto x+d$ is affine for a fixed $x$, both the functions $\underline{g_i}(x+d)-\underline{g_i}(x)$ and $\overline{g_i}(x+d)-\overline{g_i}(x)$ are convex for all $d\in \R^n$.  Hence, their pointwise maximum $\overline{\Delta g_i}(x;d)$ is also convex, $i.e,$
\begin{align}\label{g convex}
\overline{\Delta g_i}(x;d_\lambda)\le \lambda\overline{\Delta g_i}(x;d_1)+(1-\lambda)\overline{\Delta g_i}(x;d_2),\ \forall d_1,d_2\in \R^n,\  \lambda\in [0,1].
\end{align}
By combining the inequalities \eqref{f convex} and \eqref{g convex} and on simplification, we obtain
\[
\overline{h}_i(x;d_\lambda)\le\lambda \overline{h}_i(x;d_1)+(1-\lambda)\overline{h}_i(x;d_2).
\]
Now, taking maximum over $i=\overline{1,m}$, one has
\begin{align*}
    \psi(x;d_\lambda):=\max_{i=1,\dots,m} \overline{h}_i(x;d_\lambda)
    &\le \max_{i=1,\dots,m} \big \{\lambda\overline{h}_i(x;d_1)+(1-\lambda)\overline{h}_i(x;d_2)\big \}\\
    &\le \lambda \max_{i=1,\dots,m}\overline{h}_i(x;d_1)+(1-\lambda)\max_{i=1,\dots,m}\overline{h}_i(x;d_2),\\
    &= \lambda \psi(x;d_1)+(1-\lambda)\psi(x;d_2),
\end{align*}
for all $d_1, d_2\in\R^n$ and $\lambda\in[0,1]$. Therefore, $\psi(x;d)$ is a convex function.\\

\noindent Also, by the definition of $\psi(x;d)$, we have
\begin{align*}
    \psi(x;0)&=\max_{i=1,\dots,m} \overline{h}_i(x;0),\\
    &=\max_{i=1,\dots,m}\{\overline{\nabla_{gH}f_i(x)^\top\odot 0\ \oplus \Delta g_i(x;0)}\}, \\
    &=0.
\end{align*}
Hence, the proof is complete.
\end{proof}
\noindent Now, we consider the following unconstrained optimization problem at a given point $x\in dom\ F$, as follows
\begin{align*}
\textbf{OP}\quad  \min_{d\in\R^n} \phi(x;d):= \psi(x;d)+\frac{1}{2}\|d\|^2.
\end{align*}
\textbf{Note:} From Lemma~\ref{Convex lemma}, it is obvious that the objective function is strongly convex and thus, it has a unique optimal solution.\\ 

\noindent Let $\alpha(x)$ and $\beta(x)$ be the unique optimal solution and the optimal value  of OP, respectively, at a given point $x\in dom\ F$, defined by
\begin{align}
\alpha(x)&:=\arg \min_{d\in \R^n}\phi(x;d), \label{alpha}\\
\beta(x)&:=\min_{d\in\R^n}\phi(x;d)=\phi(x; \alpha(x)).
\label{beta}
\end{align}

\begin{remark}\label{remark}
    \begin{enumerate}[label=(\roman*)]
        \item If the IVF $g_i(x)
    =\mathbf{0}$ in the  MIVOP$_C$ for all $i=\overline{1,m}$, then the $\alpha(x)$ reduces to the search direction of the steepest descent method proposed by Mondal and Ghosh~\cite{mondal2025steepest}. 
        \item If $F$ is a multiobjective real valued function, that is $\overline{F_i}=\underline{F_i}$ for all $i=\overline{1,m}$ yields $\nabla f_i=\underline{\nabla_{gH}f_i}=\overline{\nabla_{gH}f_i}$, also $\overline{\Delta g_i}=\underline{\Delta g_i}=g_i(x+d)-g_i(x)$. Therefore, in this case $\alpha$ reduces to 
        \[
        \arg\min_{d\in \R^n}\ \big (\max_{i=1,\dots,m} \{\nabla f_i(x)^\top d +g_i(x+d)-g_i(x)\}+\frac{1}{2}\|d\|^2\ \big ),
        \]
        which corresponds to the search direction of the proximal gradient method for real valued multiobjective optimization problem proposed by Tanabe $et\ al.$\cite{TanabePGM}.
    \end{enumerate}
\end{remark}

\noindent To obtain the value of $\alpha(x)$, we reformulate OP as follows:
\begin{align*}
    &\min_{d\in\R^n}\big (\max_{i=1,\dots,m}\overline{h}_i(x;d)+\frac{1}{2}\|d\|^2\big )\\
    &\equiv \min_{d\in \R^n}\Big (\max_{i=1,\dots,m}\Big \{\frac{1}{2}\big (  \underline{\nabla_{gH}f_i}(x)+\overline{\nabla_{gH}f_i}(x)\big )^\top d\ +\frac{1}{2}\big( \overline{\nabla_{gH}f_i}(x)- \underline{\nabla_{gH}f_i}(x)\big )^\top|d|\\
    &\quad \hspace{2.5cm}+\frac{1}{2}\big (\underline{g_i}(x+d)-\underline{g_i}(x)+\overline{g_i}(x+d)- \overline{g_i}(x)\big )\\
    &\quad \hspace{4.9cm}+\frac{1}{2}\big|\underline{g_i}(x+d)-\underline{g_i}(x)
    -\overline{g_i}(x+d)+\overline{g_i}(x)\big|\Big \}+\frac{1}{2}\|d\|^2\Big ),
\end{align*}
equivalently,
\begin{align*}
    \textbf{COP} \qquad \min_{d,z,\omega,\tau}\  &\tau+\frac{1}{2}\|d\|^2\\
    \text{subject to}\  &\big (  \underline{\nabla_{gH}f_i}(x)+\overline{\nabla_{gH}f_i}(x)\big )^\top d\ +\big( \overline{\nabla_{gH}f_i}(x)- \underline{\nabla_{gH}f_i}(x)\big )^\top z\ + \\
    &\hspace{0.5cm} \underline{g_i}(x+d)-\underline{g_i}(x)+\overline{g_i}(x+d)-\overline{g_i}(x)+\omega_i\le 2\tau,\ \ \ \forall i=\overline{1,m},\\
    &-\omega_i\le \underline{g_i}(x+d)-\underline{g_i}(x)-\overline{g_i}(x+d)+\overline{g_i}(x)\le\omega_i,\ \forall i=\overline{1,m},\\
    & -z_j\le d_j\le z_j,\ \forall j=\overline{1,n}.
\end{align*} 
The following result characterizes the Pareto criticality of the MIVOP$_C$ in terms of $\alpha(x)$ and $\beta(x)$.
\begin{theorem}\label{Theorem 1}
For any $x\in dom\ F$, 
    let $\alpha(x)$ and $\beta(x)$ be defined by relations \eqref{alpha} and \eqref{beta}, respectively. Then the following properties are satisfied:
    \begin{enumerate}[label=(\roman*)]
    \normalfont
        \item \textit{$\beta(x)\le 0$, for all $x\in dom\ F$}
        \item \textit{If $x$ is a Pareto critical point of the MIVOP$_C$, then both $\alpha(x)$ and $\beta(x)$ vanish.} 
        \item \textit{If $x$ fails to be a Pareto critical point of the MIVOP$_C$, then $\alpha(x)\neq0$ and $\beta(x)<0$.}
        \item \textit{$\alpha(x)$ is bounded on a compact subset of $int(dom\ F)$.}
        \item \textit{Both $\alpha$ and $\beta$ are continuous mappings on $int(dom\ F)$.}
    \end{enumerate}
\end{theorem}

\begin{proof}
    (i) Since $\beta(x)$ is an optimal value of OP, then for any given $x\in dom\ F$, we have
    \[
        \beta(x)\le \psi(x;d)+\frac{\|d\|^2}{2},\ \forall d\in \R^n.
    \]
    In particular, for $d=0$, one has
    \begin{align*}
    \beta(x)&\le \psi(x;0)+\frac{\|0\|^2}{2}=0.\ (using\ Lemma~\ref{Convex lemma})
    \end{align*}
\noindent (ii) Let $x\in dom\ F$ be a Pareto critical point of the MIVOP$_C$. Then from Corollary~\ref{PCP implies no hi<0}, there exists no $d\in \R^n$ such that
\[
h_i(x;d)\prec \mathbf{0}.
\]
Equivalently, for all $d\in \R^n$, $\mathbf{0}\preceq h_i(x;d)$ and therefore, 
\begin{align*}
   &\quad \max_{i=1,\dots,m}\overline{h}_i(x;d)\ge 0,\ \forall d\in \R^n,\\
   &\Rightarrow \max_{i=1,\dots,m}\overline{h}_i(x;\alpha)+\frac{\|\alpha\|^2}{2}\ge 0,\ (for\ some\ d=\alpha)\\
   &\Rightarrow \psi(x;\alpha)+\frac{\|\alpha\|^2}{2}\ge 0,\\
   &\Rightarrow \beta(x)\ge 0,
\end{align*}
    which in view of the results of part~(i) yields $\beta(x)=0$ and $\alpha(x)=0$.

\medskip
 \noindent (iii) Let $x\in dom\ F$ not be a Pareto critical point of the MIVOP$_C$, then from Proposition~\ref{NPC implies hi}, there exists $ d\in \R^n$, such that 
 \begin{align*}
 h_i&(x; d)\prec \mathbf{0}, \forall i=\overline{1,m},\\
 \text{or, }\big [\underline{h_i}&(x; d),\ \overline{h}_i(x;d)\big ]\prec \mathbf{0},\\
     \Rightarrow &\max_{i=1,\dots,m}\overline{h}_i(x;d)<0,\\
     \Rightarrow &\psi(x; d)<0.
 \end{align*}
Also, since $\psi(x;0)=0$, it follows from the above inequality $ d\neq 0$. 

\medskip
\noindent Further, since $\psi(x;d)$ is convex, for any $t\in(0,1)$, one has
\begin{align*}
    \psi(x;t d)&=\psi(x;t d+(1-t)\ 0)\\
    &\le t\psi(x; d).
\end{align*}
Therefore,
\begin{align*}
    \phi(x;t d)&=\psi(x; td)+\frac{1}{2}t^2\| d\|^2,\\
     &\le t\psi(x; d)+\frac{1}{2}t^2\| d\|^2,\\
     &= t\big (\psi(x; d)+\frac{1}{2}t\| d\|^2\big ).
\end{align*}
 Since $\psi(x; d)<0$, there exists $\overline t=\min\big \{1,-\frac{\psi(x; d)}{C}\big \}$, where $C=\frac{1}{2}\| d\|^2$, such that for every $t\in (0,\overline t)$, $\psi(x; d)+\frac{1}{2}t\| d\|^2<0$, and we have
 $$\phi(x;t d)<0.$$
 Now, since $\beta(x)$ is the optimal value of OP, one has
 \begin{align*}
     \beta(x)\le \phi(x;td)<0.
 \end{align*}
 Finally, since $\phi(x;0)=0\neq\beta(x)$, thus the point $0\in \R^n$ does not minimize the objective function $\phi(x;d)$, for a given point $x$, and we conclude that $\alpha(x)\neq 0$.

\medskip
\noindent (iv) Let $\mathcal T\subset int(dom\ F)$ be a compact set. Since the IVF $\nabla_{gH}f_i$ is $gH$-Lipschitz continuous, from Lemma~\ref{lemma gh lip}, its endpoint functions are also Lipschitz continuous, and therefore, $\underline{\nabla_{gH}f_i}$ and $\overline{\nabla_{gH}f_i}$ are bounded on the compact set $\mathcal T$, for all $i=\overline{1,m}$. Consequently, there exists a constant $J>0$, such that
 \begin{align*}
 \max\{\|\underline{\nabla_{gH}f_i}(x)\|,\|\overline{\nabla_{gH}f_i}(x)\|\}\le J,\ \forall x\in \mathcal T,\ i=\overline{1,m}.
 \end{align*}
Now, using the Cauchy-Schwarz inequality, for a given $x\in\mathcal T$, one has
\begin{align*}
    \frac{1}{2}\big (  \underline{\nabla_{gH}f_i}(x)+\overline{\nabla_{gH}f_i}(x)\big )^\top d&\ge -\Big |\frac{1}{2}\big (  \underline{\nabla_{gH}f_i}(x)+\overline{\nabla_{gH}f_i}(x)\big )^\top d\Big |,\\
    &\ge -\frac{1}{2}\big (\|\underline{\nabla_{gH}f_i}(x)\|+\|\overline{\nabla_{gH}f_i}(x)\|\big )\|d\|,\\
    &\ge -J\|d\|,\ \forall d\in \R^n,\ i=\overline{1,m}.
\end{align*}
Similarly, 
\begin{align*}
    \frac{1}{2}\big( \overline{\nabla_{gH}f_i}(x)- \underline{\nabla_{gH}f_i}(x)\big )^\top|d|\ge -J\|d\|,\  \forall i=\overline{1,m}.
\end{align*}
On adding the above two inequalities, we obtain
\begin{align}
    \frac{1}{2}\big (  \underline{\nabla_{gH}f_i}(x)+\overline{\nabla_{gH}f_i}(x)\big )^\top d\ +\frac{1}{2}\big( \overline{\nabla_{gH}f_i}(x)- \underline{\nabla_{gH}f_i}(x)\big )^\top|d|\ge -2J\|d\|.
    \label{CS ine fi}
\end{align}
 Further, since the endpoint functions $\underline{g_i}$ and $\overline{g_i}$ of the IVF $g_i$ belong to $\Gamma_0(\R^n)$, the functions $\underline{g_i}$ and $\overline{g_i}$ have bounded subgradients on the compact set $\mathcal T$, for all $i=\overline{1,m}$ (see, for instance Theorem~24.7 of Rockafellar~\cite{tyrrell1970convex}) and thus, there exists a constant $K>0$ such that
\[
\max\big \{\|\underline{\xi_i}\|,\|\overline{\xi_i}\|\big \} \le K,
\]
 where $\underline{\xi_i}$ and $\overline{\xi_i}$ are the subgradients of $\underline{g_i}$ and $\overline{g_i}$ at the point $x\in \mathcal T$, respectively, for each $i=\overline{1,m}$ and therefore, for every $d\in \R^n$, we have
\begin{align*}
    \underline{g_i}(x+d)-\underline{g_i}(x)\ge\langle \underline{\xi_i},d\rangle\
  \text{and }  \overline{g_i}(x+d)-\overline{g_i}(x)\ge\langle \overline{\xi_i},d\rangle.
\end{align*}
Again, using the Cauchy-Schwarz inequality and proceeding as above, we have
\begin{align*}
    \underline{g_i}(x+d)-\underline{g_i}(x)\ge-K\|d\|\ \text{and }
    \underline{g_i}(x+d)-\underline{g_i}(x)\ge-K\|d\|,\ \forall d\in \R^n,\ i=\overline{1,m}. 
\end{align*}
Now, since $\overline{\Delta g_i}(x;d)=\max\{\underline{g_i}(x+d)-\underline{g_i}(x), \underline{g_i}(x+d)-\underline{g_i}(x)\}$, we obtain 
 \begin{align}
 \overline{\Delta g_i}(x;d)\ge -K\|d\|, \ \forall i=\overline{1,m}.
 \label{CS gi}
\end{align}
 On combining the inequalities \eqref{CS ine fi} and \eqref{CS gi}, we obtain
 \[
 \overline{h}_i(x;d)\ge-2J\|d\|-K\|d\|=-\mathcal C\|d\|,\ \forall d\in \R^n,i=\overline{1,m},
 \]
where $\mathcal C=2J+K$, which in the view of the expression~\eqref{psi def} yields
\begin{align*}
\psi(x;d)&\ge -\mathcal C\|d\|,\ \forall d\in\R^n,\\
\text{or, }  \psi(x;d)+\frac{1}{2}\|d\|^2&\ge -\mathcal C\|d\|+\frac{1}{2}\|d\|^2,\\
&\ge \min\big \{ -\mathcal C\|d\|+\frac{1}{2}\|d\|^2\big \},\\
&=-\frac{\mathcal {C}^2}{2}=:-\mathcal C_0,
\end{align*}
which upon taking the minimum over $d\in \R^n$ gives
\begin{align}
 -\beta(x)\le \mathcal C_0.
\label{beta C0}
\end{align}
Now, since $\phi(x;d)$ is strongly convex and $\alpha(x)$ is its unique minimizer at a given point $x\in \mathcal T$, then for all $\lambda\in (0,1]$ and $d\in \R^n$, we have
\begin{align*}
    \phi(x;\alpha(x))&\le \phi(x;\lambda d+(1-\lambda)\alpha(x)),\\
    &\le \lambda \phi(x;d)+(1-\lambda)\phi(x;\alpha(x))-\frac{1}{2}\lambda(1-\lambda)\|d-\alpha(x)\|^2.
\end{align*}
Rearranging the above inequality and taking the limit $\lambda\to0^+$, we obtain
\begin{align*}
\phi(x;d)&\ge \phi(x;\alpha(x))+\frac{1}{2}\|d-\alpha(x)\|^2.
\end{align*}
Since $\phi(x;0)=0$, we have
\begin{align*}
0&\ge \phi(x;\alpha(x))+\frac{1}{2}\|\alpha(x)\|^2,\\
&=\beta(x)+\frac{1}{2}\|\alpha(x)\|^2, \\
   \text{or,}\quad   \frac{1}{2}\|\alpha(x)\|^2&\le -\beta(x)\le \mathcal C_0,\ (from\ inequality\ \eqref{beta C0}) \\
   \text{or,}\quad    \|\alpha(x)\|&\le \sqrt{2\mathcal C_0}=:\mathcal K,\quad \forall x\in \mathcal T.
 \end{align*}
Hence, $\alpha(x)$ is bounded on the compact set $\mathcal T$.

\medskip
\noindent (v) Since the IVF $f_i$ has $gH$-Lipschitz continuous gradients, for each $i=\overline{1,m}$, the gradients $ \underline{\nabla_{gH}}f_i(x)$ and $\overline{\nabla_{gH}}f_i(x)$ are continuous for all $x\in int(dom\ F)$. Further, the map $d\mapsto (d,|d|)$ is continuous for all $d\in \R^n$. Hence,
\[
\frac{1}{2}\big (\underline{\nabla_{gH}}f_i(x)+\overline{\nabla_{gH}}f_i(x)\big )^\top d+\frac{1}{2}\big (\overline{\nabla_{gH}}f_i(x)-\underline{\nabla_{gH}}f_i(x)\big)|d|,
\]
is continuous in $(x,d)\in int(dom\ F)\times \R^n$.

\medskip
\noindent Also, from Lemma~\ref{gi gh-continuous}, the endpoint functions $\underline{g}_i,\overline{g}_i$ of the IVF $g_i$ are continuous, thus $\overline{\Delta g_i}(x;d)$ is continuous in $(x,d)\in int(dom\ F)\times \R^n$ and therefore, it follows that the objective function of the OP, $i.e,$
\[
\phi(x;d)=\max_{i=1,\dots,m}\overline{h}_i(x;d)+\frac{1}{2}\|d\|^2,
\]
 is continuous in $(x,d)\in int(dom\ F)\times \R^n$, where
 \[
\overline{h}_i(x;d)=\frac{1}{2}\big (\underline{\nabla_{gH}}f_i(x)+\overline{\nabla_{gH}}f_i(x)\big )^\top d+\frac{1}{2}\big (\overline{\nabla_{gH}}f_i(x)-\underline{\nabla_{gH}}f_i(x))|d|+\overline{\Delta g_i}(x;d),\ \forall i=\overline{1,m}.
\]
\noindent Now, consider an arbitrary sequence $\{x^k\}$ converging to $x^*\in int(dom\ F)$, then there exists $k_0\in\mathbb{N}$ such that
$x^k\in B(x^*;r)$, $r>0$, for all $k\ge k_0$, where $B(x^*;r)=\{x\in\mathbb{R}^n:\|x-x^*\|\le r\}\subset int(dom\ F)$ and thus, from the conclusion of part~(iv) it follows that the sequence $\{\alpha(x^k)\}$ is bounded for all $k\ge k_0$, and therefore by Bolzano-Wierstrass Theorem, there exist a subsequence
$\{\alpha(x^{{k_j}})\}$ converging to $\overline \alpha\in \R^n$, as $j\to \infty$. Also, since $\phi(x;d)$ is continuous in $(x,d)\in int(dom\ F)\times \R^n$, then the following implications hold:
\begin{equation}
\begin{aligned}
\phi(x^{k_j};\alpha(x^{k_j}))&\to \phi(x^*;\overline \alpha) \\
\text{and } \phi(x^{k_j};\alpha(x^*))&\to \phi(x^*;\alpha(x^*))=\beta(x^*).
\label{phi star 1}
\end{aligned}
\end{equation}
On the other hand, since $\alpha(x^{k_j})$ is the minimizer of $\phi(x^{k_j};d)$ for all $d\in \R^n$, we have
\begin{align*}
\phi(x^{k_j};\alpha(x^{k_j})) &\le \phi(x^{k_j};d),\\
   \Rightarrow \phi(x^*;\overline \alpha)&\le \phi(x^*;d),\ \text{as }j\to \infty,\ \forall d\in \R^n.
\end{align*}
Consequently, $\overline{\alpha}$ is the minimizer of the function $\phi(x^*;d)$. Therefore,
\begin{align}
   \phi(x^*;\overline \alpha)=\min_{d\in\R^n}\phi(x^*;d):=\beta(x^*).
   \label{phi eq 2}
\end{align}
Since, the function $\phi(x^*;d)$ has a unique minimizer $\alpha(x^*)$, it follows from the relations \eqref{phi star 1} and \eqref{phi eq 2} $\overline \alpha=\alpha(x^*)$. Thus, every convergent subsequence of the bounded sequence $\{\alpha(x^k)\}$ converges to $\alpha(x^*)$, and we conclude that
\begin{align*}
   \alpha(x^k)\to \alpha(x^*).
\end{align*}
Therefore,
\begin{align*}
\lim_{k\to\infty}\beta(x^k)&=\lim_{k\to\infty}\phi(x^k;\alpha(x^k)),\\
    &=\phi(x^*;\alpha(x^*)):=\beta(x^*).
\end{align*}
From the above two sequential criterion, we conlude that the mappings $x\mapsto \alpha(x)$ and $x\mapsto \beta(x)$ are continuous on $int(dom\ F)$.
\end{proof}

\begin{remark}\label{theorem 3.1 remark}
    In the above Theorem~\ref{Theorem 1}, assertions (ii) and (iii) provides complementary characterizations of Pareto criticality of the MIVOP$_C$. In particular, (iii) is the contrapositive of (ii), and hence the two statements are equivalent. Therefore, a point $x\in \R^n$ is a Pareto critical point of MIVOP$_C$ if and only if $\alpha(x)=0$ and $\beta(x)=0$.
\end{remark}

\noindent We now propose the proximal gradient method for locating the Pareto critical points of the MIVOP$_C$. This method generates a sequence $\{x^k\}_{k\in \mathbb{N}}$ under the iterative procedure given as follows:
\[
x^{k+1}=x^k+t_k d^k,
\]
where $t_k$ is the step size and $d^k$ is the search direction. In particular, the optimal solution $\alpha(x^k)$ obtained by solving the OP is the search direction $d^k$, $i.e,$ $d^k=\alpha(x^k)$ at the given point $x=x^k$, which is generally called the proximal gradient direction at iteration $k$. The proposed method is presented below:
\begin{algorithm}[H]
\caption{Proximal gradient method for solving MIVOP$_C$}
\label{alg:1}
\begin{algorithmic}[1]
\State \textbf{Input:} Provide $F=(F_1,F_2,\dots,F_m)^\top$ with each of its components given as $F_i=f_i\oplus g_i$ for all $i=\overline{1,m}$, where $f_1,f_2,\dots,f_m$ are $gH$-differentiable IVFs with $gH$-Lipschitz continuous gradients and $g_1,g_2,\dots,g_m$ are proper convex and $gH$-lower semicontinuous IVFs.
\State Choose the initial point $x^0\in dom\ F$, predefined parameter $\sigma\in(0,1)$, the contraction factor $\gamma\in(0,1)$, the stopping tolerance $\epsilon>0$. Initialize $k=0$.

\State At the current iterate $x^k$, evaluate the $gH$-gradient of the smooth part $f_i$ of each component $F_i$, namely,
\[
\nabla_{gH}f_i(x^k)=\left [\underline{\nabla_{gH}f_i}(x^k),\overline{\nabla_{gH}f_i}(x^k)\right],\ \forall i=\overline{1,m}
.\]

\State Compute the proximal gradient direction $d^k$ and optimal value $\beta^k$ of the subproblem OP at $x=x^k$.

\State If $\beta^k\ge-\epsilon$, terminate the procedure and designate $x^k$ as an approximate Pareto critical point. Otherwise, go to Step 6.

\State Set $t_k=1$. Reduce the step size $t_k:=\gamma t_k$ until the Armijo-like condition
\begin{align}\label{armijo 1}
F_i(x^k+t_kd^k)\preceq F_i(x^k)\oplus [\sigma t_k,\sigma t_k]\odot \overline{h}_i(x^k;d^k)
\end{align}
is satisfied.
\State Update the iterate by 
\[
x^{k+1}=x^{k}+t_kd^k.
\]
Set $k=k+1$ and return to \textbf{Step~3}
\end{algorithmic}
\end{algorithm}

\noindent Note that the line search condition \eqref{armijo 1} in the above Algorithm can be equivalently written as:
\begin{equation}\label{armijo 2}
\begin{aligned}
    \underline{F_i}(x^k+t_kd^k)&\le \underline{F_i}(x^k)+\sigma t_k\ \overline{h_i}(x^k;d^k) \\
   \text{and }\quad \overline{F_i}(x^k+t_kd^k)&\le \overline{F_i}(x^k)+\sigma t_k\ \overline{h_i}(x^k;d^k)
\end{aligned}
\end{equation}
\noindent The following results guarantees that the line search in Algorithm~\ref{alg:1} is well defined.

\begin{theorem}\label{t exist}
Let $d^k$ be the proximal gradient direction obtained by solving the OP in Algorithm~\ref{alg:1}, computed at $x=x^k$. If $x^k$ is not a Pareto critical point of the MIVOP$_C$, then there exists $\overline t_k>0$ such that 
\begin{align*}
F_i(x^k+td^k)\prec F_i(x^k)\oplus [\sigma t,\sigma t]&\odot \overline{h_i}(x^k;d^k),
\end{align*}
for all $t\in(0,\overline t_k\ ]$ and $\sigma\in(0,1)$.
\end{theorem}
\begin{proof}
Let $x^k\in dom\ F$ be not a Pareto critical point of the MIVOP$_C$, then from the third assertion of Theorem~\ref{Theorem 1}, we obtain $ \beta(x^k)<0$ and thus, it follows that
\begin{align*}
     \psi(x^k;d^k)+\frac{\|d^k\|^2}{2}&<0
    \Rightarrow \psi(x^k;d^k)<0
    \Rightarrow h_i(x^k;d^k)\prec\mathbf{0}.\ (from\ equation~\eqref{implies})
\end{align*}
Further, since the IVF $F_i$ admits a composite structure for all $i=\overline{1,m}$, therefore, for any $t>0$, we have
\begin{align*}
    \lim_{t\to 0^+}&\frac{1}{t}\odot\big (F_i(x^k+td^k)\ominus_{gH}F_i(x^k)\big )\\
    &=  \lim_{t\to 0^+}\frac{1}{t}\odot \Big (\big (f_i(x^k+td^k)\oplus g_i(x^k+td^k)\big)\ominus_{gH}\big(f_i(x^k)\oplus g_i(x^k)\big)\Big ),\\
    &=\lim_{t\to 0^+}\frac{1}{t}\odot \Big (\big (f_i(x^k+td^k)\ominus_{gH} f_i(x^k)\big)\oplus\big(g_i(x^k+td^k)\ominus_{gH} g_i(x^k)\big)\Big ),\\
    &=\lim_{t\to 0^+}\frac{1}{t}\odot\big (f_i(x^k+td^k)\ominus_{gH} f_i(x^k)\big )\oplus\lim_{t\to 0^+}\frac{1}{t}\odot\big (g_i(x^k+td^k)\ominus_{gH} g_i(x^k)\big ),\\
        &\preceq \nabla_{gH}f_i(x^k)^\top\odot d^k\ \oplus\ \big (g_i(x^k+d^k)\ominus_{gH}g_i(x^k)\big),\ (using\ Lemmas~\ref{gh diff lemma}\ and~\ref{DD Lemma}) \\
        &=h_i(x^k;d^k)\prec \mathbf{\sigma}\ \odot h_i(x^k;d^k).\quad (since\  h_i(x^k;d^k)\prec\mathbf{0}\ and \ \sigma\in(0,1))
\end{align*}
Consequently, there exist $\overline t_k>0$ such that for all $t\in (0,\overline t_k \ ] $ and $i=\overline{1,m}$, we have 
\begin{align*}
 F_i(x^k+td^k)\ominus_{gH} F_i(x^k)&\prec \sigma t\ \odot h_i(x^k;d^k),\\
\text{or, } F_i(x^k+td^k) &\prec F_i(x^k)\oplus \sigma t\ \odot h_i(x^k;d^k),\\
&\preceq F_i(x^k)\oplus \sigma t \odot \big[\ \overline{h_i}(x^k;d^k),\overline{h_i}(x^k;d^k)\big ],\\
\text{or, } F_i(x^k+td^k) &\prec F_i(x^k)\oplus [\sigma t,\sigma t]\ \odot \overline{h_i}(x^k;d^k),
\end{align*}
which concludes the proof.
\end{proof}

\begin{lemma}\label{lemma f lipsch}
    Let $f_i:\R^n\to \IR$ be a $gH$-differentiable IVF with $gH$-Lipschitz continuous gradients having Lipschitz constant $L_i>0$, for all $i=\overline{1,m}$. Then the following inequalities holds for every $x,d\in \R^n$,
    \begin{equation}\label{fi Lipsch}
    \left.
    \begin{aligned}
      \underline{f_i}(x+d)-\underline{f_i}(x)&\le \frac{L_i}{2}\|d\|^2+\overline{\nabla_{gH}f_i}(x)^\top d\\
      \text{and }\ \overline{f_i}(x+d)-\overline{f_i}(x)&\le \frac{L_i}{2}\|d\|^2+\overline{\nabla_{gH}f_i}(x)^\top d. 
    \end{aligned}
    \right.
    \end{equation}
\end{lemma}
\begin{proof}
The proof is analogous to the proof of Lemma~3.3 of Mondal and Ghosh~\cite{mondal2025steepest} and hence been omitted. 
\end{proof}
\noindent To establish a positive lower bound for the step sizes  generated by Algorithm~\ref{alg:1}, we present the following result.

\begin{theorem}\label{t minimum}
Let $x^k\in dom\ F$ be a non Pareto critical iterate, and let
    $d^k$ denote the corresponding proximal gradient direction obtained by solving the OP. Suppose that $\overline t_k$ is the step length satisfying the Armijo-like condition~\eqref{armijo 2} in Algorithm~\ref{alg:1}, $\overline t_k$ admits the following uniform estimate:
    \begin{align}
    \overline t_k\ge t_{\min}:=\min \left \{\frac{\gamma(1-\sigma)}{L_{\max}},1\right \},
    \label{t minimum equation}
    \end{align}
    where $L_{\max}=\max\limits_{i=1,\dots,m}L_i$ and $\sigma,\gamma\in(0,1)$.
\end{theorem}

\begin{proof}
    Since $f_i:\R^n\to \IR$ is $gH$-differentiable IVF with $gH$-Lipschitz continuous gradients having $L_i>0$ Lipschitz constants for all $i=\overline{1,m}$, then for any $x^k,d^k\in \R^n$, Lemma~\ref{lemma f lipsch} yields 
    \begin{equation}\label{fi Lipsch}
    \left.
    \begin{aligned}
      \underline{f_i}(x^k+td^k)-\underline{f_i}(x^k)&\le \frac{L_it^2}{2}\|d^k\|^2+t\ \overline{\nabla_{gH}f_i}(x^k)^\top d^k\\
    \text{and }  \overline{f_i}(x^k+td^k)-\overline{f_i}(x^k)&\le \frac{L_it^2}{2}\|d^k\|^2+t\ \overline{\nabla_{gH}f_i}(x^k)^\top d^k,\ \forall t\in(0,1]
    \end{aligned}
    \right.
    \end{equation}
 Further, since $g_i,\ i=\overline{1,m}$ is a convex IVF, from Lemma~\ref{convex IVF lemma}, the endpoint functions $\underline{g_i}$ and $\overline{g_i}$ are also convex and therefore
    \begin{align*}
        \underline{g_i}(x^k+td^k)&=\underline{g_i}(t(x^k+d^k)+(1-t)x^k)\\
        &\le t\underline{g_i}(x^k+d^k)+(1-t)\underline{g_i}(x^k)\\
        &=\underline{g_i}(x^k)+t\big (\underline{g_i}(x^k+d^k)-\underline{g_i}(x^k)\big ),  \forall i=\overline{1,m},\  t\in(0,1].
    \end{align*}
    Similarly, the upper endpoint function $\overline{g_i}$ admits the same inequality and therefore
    \begin{equation*}
    \begin{aligned}
   \overline{g_i}(x^k+td^k)-\overline{g_i}(x^k)&\le t\big (\overline{g_i}(x^k+d^k)-\overline{g_i}(x^k)\big ),\ \forall i=\overline{1,m},\ t\in(0,1]
    \end{aligned}
    \end{equation*}
On setting $ \overline{\Delta g_i}(x^k;d^k)=\max\{\underline{g_i}(x^k+d^k)-\underline{g_i}(x^k),\ \overline{g_i}(x^k+d^k)-\overline{g_i}(x^k)\}$ for all $i=\overline{1,m}$, it follows from the above two inequalities
\begin{equation}\label{eq gi}
\begin{aligned}
 \underline{g_i}(x^k+td^k)-\underline{g_i}(x^k)&\le t\ \overline{\Delta g_i}(x^k:d^k)\\
\text{and } \overline{g_i}(x^k+td^k)-\overline{g_i}(x^k)&\le t\ \overline{\Delta g_i}(x^k;d^k).
\end{aligned}
\end{equation}
\noindent Now, on combining the inequalities~\eqref{fi Lipsch} and \eqref{eq gi}, we obtain
\begin{equation}\label{eq 8.}
\begin{aligned}
 \underline{F_i}(x^k+td^k)-\underline{F_i}(x^k)&\le \frac{L_it^2}{2}\|d^k\|^2+t\ \overline{h_i}(x^k;d^k)\\
\text{and } \overline{F_i}(x^k+td^k)-\overline{F_i}(x^k)&\le \frac{L_it^2}{2}\|d^k\|^2+t\ \overline{h_i}(x^k;d^k),\ \forall i=\overline{1,m},\ \forall t\in(0,1].
\end{aligned}
\end{equation}

\noindent On the other hand, since $\overline t_k$ satisfies the Armijo-like conditions~\eqref{armijo 2}, it is obvious that $\dfrac{\overline t_k}{\gamma}$ does not satisfy the conditions~\eqref{armijo 2}, and therefore

\begin{equation}\label{eq 9}
\begin{aligned}
    \underline{F_{i}}\left (x^k+\frac{\overline t_k}{\gamma}d_k\right )- \underline{F_{i}}(x^k)&> \sigma \frac{\overline t_k}{\gamma}\overline{h_{i}}(x^k;d^k)\\
    \text{and }\overline{F_{i}}\left (x^k+\frac{\overline t_k}{\gamma}d^k\right )-\overline{F_{i}}(x^k)&> \sigma \frac{\overline t_k}{\gamma}\overline{h_{i}}(x^k;d^k),\ \forall i=\overline{1,m}.
\end{aligned}
\end{equation}

\noindent Since $\gamma\in (0,1)$ and $\overline{t}_k$ takes the value $1,\gamma,\gamma^2,\dots$, we have $\dfrac{\overline t_k}{\gamma}\in(0,1]$. Substituting $t=\dfrac{\overline t_k}{\gamma}$ in \eqref{eq 8.} and using the above inequalities, we obtain
\begin{align*}
    &\sigma\frac{\overline t_k}{\gamma}\overline{h_{i}}(x^k;d^k)<\frac{L_{i}\overline t_k^2}{2\gamma^2}\|d^k\|^2+\frac{\overline t_k}{\gamma}\overline{h}_{i}(x^k;d^
    k),\\
    \text{or, } &0<\frac{L_{i} \overline t_k}{2\gamma}\|d^k\|^2+(1-\sigma)\overline{h_{i}}(x^k;d^k),\ \forall i=\overline{1,m}.
\end{align*}
Now, on taking maximum over $i=\overline{1,m}$, we obtain
\begin{align*}
     0&<\frac{L_{\max}\ \overline t_k}{2\gamma}\|d^k\|^2+(1-\sigma)\max_{i=1,\dots,m}\overline{h_i}(x^k;d^k),\\
     &=\frac{L_{\max}\ \overline t_k}{2\gamma}\|d^k\|^2+(1-\sigma)\psi(x^k;d^k),\\
     &<\frac{L_{\max}\ \overline t_k}{2\gamma}\|d^k\|^2+(1-\sigma)\big (-\frac{1}{2}\|d^k\|^2\big ),\ (from\  \text{Theorem}~\ref{Theorem 1}(iii)) \\
   \text{or, } \ &\overline t_k>\frac{\gamma (1-\sigma)}{L_{\max}}.
\end{align*}
Since $\overline t_k\le 1$, it follows from the above inequality
\begin{align*}
   \overline t_k\ge t_{\min}:=\min \left \{\frac{\gamma(1-\sigma)}{L_{\max}},1\right \}.
\end{align*}
Hence, the proof is complete.
\end{proof}

\section{Convergence analysis}\label{section 4}

\noindent In this section, we discuss about the convergence analysis of the proposed Algorithm~\ref{alg:1} and derive that the sequence generated by the Algorithm~\ref{alg:1} converges to a Pareto critical point of the MIVOP$_C$  followed by the rate of convergence in the nonconvex setting.

\medskip
\noindent The following result given in Chen and Teboulle~\cite{Chen proximal} forms the fundamental tool for studying the convergence of proximal type algorithms.

\begin{theorem}[\cite{Chen proximal}]
\label{theorem tpp}
    Let $\Theta:\R^n\to\R_{\infty}$ be a proper convex function and let $\hat d\in \R^n$ such that 
    \[
    \hat d=\arg \min_{d\in\R^n}\big \{\Theta(d)+\frac{1}{2}\|d-e\|^2\big \},\ \forall e\in \R^n.
    \]
    Then for all $l\in \R^n$, we have
    \begin{align*}
    \Theta(\hat d)-\Theta(l)\le -\frac{1}{2}\|l-\hat d\|^2-\frac{1}{2}\|e-\hat d\|^2+\frac{1}{2}\|l-e\|^2.
    \end{align*}
\end{theorem}

\begin{lemma}\label{Lemma 4.1}
    Let $\{d^k\}$ be the sequence of proximal gradient directions generated by Algorithm~\ref{alg:1}, computed at the iterates $x^k\in dom\ F$, then 
    \[
    \psi(x^k;d^k)\le -\|d^k\|^2,\quad \forall k\ge 0.
    \]
\end{lemma}

\begin{proof}
   Define $\Theta(d):=\psi(x^k;d)$ for a given iterate $x^k\in dom\ F$. Now, since $d^k$ is the optimal solution of OP at a given iterate $x^k$, we have
    \begin{align*}
    d^k&=\arg \min_{d\in\R^n}\left\{\psi(x^k;d)+\frac{1}{2}\|d\|^2\right\},\\
     &=\arg \min_{d\in\R^n}\left\{\Theta(d)+\frac{1}{2}\|d-0\|^2\right\},
    \end{align*}
     which by using Theorem~\ref{theorem tpp}, yields
   \begin{align*}
       \Theta(d^k)-\Theta(0)&\le -\frac{1}{2}\|0-d^k\|^2-\frac{1}{2}\|0-d^k\|^2+\frac{1}{2}\|0-0\|^2,\ (for\ some\ 0\in\R^n)\\
       &= -\|d^k\|^2.
\end{align*}
       Also, since $\Theta(0):= \psi(x^k;0)=0$, it follows from the above inequality
       \begin{align*}
       \psi(x^k;d^k)&\le -\|d^k\|^2,\ \forall k\ge 0.
       \qedhere
   \end{align*}
\end{proof}
\noindent We now establish the convergence of the sequence of the proximal gradient directions.

\begin{theorem}\label{tk propositon}
 Let $\{d^k\}$ be the sequence of proximal gradient directions generated by Algorithm~\ref{alg:1}, computed at the iterates $x^k\in dom\ F$. Then the following assertions hold:
 \begin{enumerate}[label=(\roman*)]
\item $F_i(x^{k+1})\preceq F_i(x^k)$,  for every $i=\overline{1,m}$ and $k\ge 0$.
\item If the sequence $\{F_i(x^k)\}$ is bounded below for all $i=\overline{1,m}$, then
   \[
    \sum_{k=0}^\infty\|d^k\|^2<\infty,
    \]
    where $t_k$ is the step size satisfying the Armijo-like condition~\eqref{armijo 2}.
 \end{enumerate}
\end{theorem}

\begin{proof}
 (i) Since $t_k$ satisfies the Armijo-like condition~\eqref{armijo 2}, we have
    \begin{equation*}
    \begin{aligned}
        \underline{F_i}(x^k+t_kd^k)&\le \underline{F_i}(x^k)+\sigma t_k\overline{h_i}(x^k;d^k)\\
        \text{and }\ \overline{F_i}(x^k+t_kd^k)&\le \overline{F_i}(x^k)+\sigma t_k\overline{h_i}(x^k;d^k),\ \forall i=\overline{1,m},
   \end{aligned}
    \end{equation*}   
which by  Lemma~\ref{Lemma 4.1} and the  fact $\overline{h_i}(x^k;d^k)\le \psi(x^k;d^k)$ for all $i=\overline{1,m}$, gives
\begin{equation}
   \begin{aligned}
     \underline{F_i}(x^k+t_kd^k)&\le \underline{F_i}(x^k)-\sigma t_k\|d^k\|^2 \\ 
     \text{and }\ \overline{F_i}(x^k+t_kd^k)&\le \overline{F_i}(x^k)-\sigma t_k \|d^k\|^2,\ \forall i=\overline{1,m}
     \label{ineq eq}
\end{aligned}
\end{equation}
Also, since $\sigma\in(0,1)$ and $t_k>0$, it follows from the above inequality
  \begin{equation*}
   \begin{aligned}   
     \underline{F_i}(x^{k+1})\le \underline{F_i}(x^k)&\text{ and }\overline{F_i}(x^{k+1})\le \overline{F_i}(x^k),\ \forall i=\overline{1,m},
     \end{aligned}
\end{equation*}
equivalently,
\begin{align*}
     F_i(x^{k+1})&\preceq F_i(x^k), \ \forall i=\overline{1,m}\text{ and } k\ge0.
     \end{align*}
     
\noindent (ii) Let the sequence $\{F_i(x^k)\}$ be bounded from below, then there exists a lower bound $\mathcal F_i:=[\ \underline{\mathcal F_i}\, , \overline{F_i}\ ]\in \IR$ such that $\mathcal F_i\preceq F_i(x^k)$ for all $i$ and $k$, and thus on summing up the inequalities~\eqref{ineq eq} from $k=0$ to $k=\hat k$, where $\hat k\in \mathbb{N}$, we obtain
\begin{align*}
        \underline{\mathcal F_i}\le\underline{F_i}(x^{\hat k+1})&\le \underline{F_i}(x^0)-\sigma\sum_{k=0}^{\hat k} t_k\|d^k\|^2\\
      \text{and } \ \overline{\mathcal F_i}\le \overline{F_i}(x^{\hat k+1})&\le \overline{F_i}(x^0)-\sigma\sum_{k=0}^{\hat k} t_k\|d^k\|^2,\ \forall i=\overline{1,m}.
\end{align*}
Therefore,
\begin{align*}
   \sum_{k=0}^{\hat k} t_k\|d^k\|^2\le \frac{\underline{F_i}(x^0)-\underline{\mathcal F_i}}{\sigma}\ \ 
\text{and }\  \sum_{k=0}^{\hat k} t_k\|d^k\|^2\le \frac{\overline{F_i}(x^0)-\overline{\mathcal F_i}}{\sigma},\ \forall i=\overline{1,m},
\end{align*}
which by taking $\hat k\rightarrow \infty$, yields
\[
 \sum_{k=0}^{\infty} t_k\|d^k\|^2<\infty.
\]
Further, since $t_k$ satisfies the Armijo-like condition~\eqref{armijo 2}, then by Theorem~\ref{t minimum}, there exists a uniform lower bound $t_{\min}$ on the step length $t_k$ such that $t_k\ge t_{\min}$ for all $k$, and therefore it follows from the above inequality
 \begin{align*}
    t_{\min}\sum_{k=0}^{\infty}&\|d^k\|^2\le \sum_{k=0}^{\infty} t_k\|d^k\|^2<\infty,\\[4pt]
    \Rightarrow\ &\sum_{k=0}^{\infty}\|d^k\|^2<\infty.
    \qedhere
    \end{align*}
\end{proof}

\noindent In the following result, we establish the global convergence of the sequence $\{x^k\}$ generated by the Algorithm~\ref{alg:1} to a Pareto critical point of MIVOP$_C$.

\begin{theorem}\label{thm:global}
Let $\{x^k\}$ be the sequence generated by Algorithm~\ref{alg:1}. Then every
accumulation point of $\{x^k\}$ lying in $int(dom\ F)$ is a Pareto critical point of the MIVOP$_C$.
\end{theorem}

\begin{proof}
If $\epsilon=0$ and the Algorithm~\ref{alg:1} stops at some iteration $k$, then $\beta(x^k)\ge 0$. However,  from the first assertion of Theorem~\ref{Theorem 1}, we have $\beta(x^k)\le0$. Consequently,
$\beta(x^k)=0$ which yields $x^k$ is a
Pareto critical point of MIVOP$_C$ by the third assertion of Theorem~\ref{Theorem 1}. So, we assume that Algorithm~\ref{alg:1} generates an infinite sequence $\{x^k\}$ with
$\beta(x^k)<0$ for all $k$.

\medskip
\noindent Let $\overline{x}\in int(dom\ F)$ be an arbitrary accumulation point of the sequence $\{x^k\}$, then there exists a subsequence $\{x^{k_j}\}$ such that
$x^{k_j}\to\overline{x}$ as $j\to\infty$. Further, since the IVF $f_i$ is $gH$-differentiable on $\R^n$ and the IVF $g_i$ is proper convex, therefore from Lemma~\ref{gi gh-continuous}, $g_i$ is $gH$-continuous on $int(dom\ F)\subseteq \R^n$ for all $i=\overline{1,m}$. Consequently, the IVF 
${F}_i={f_i}\oplus{g}_i$ is $gH$-continuous at $\overline{x}\in int(dom\ F)$, and therefore
\[
{F}_i(x^{k_j})\to{F}_i(\overline{x})\quad\text{as } j\to\infty,
\quad \forall i=\overline{1,m}.
\]
Now, let $k\in\mathbb{N}$ be arbitrary such that $k_j\ge k$ for some $j\in\mathbb{N}$, then by the first assertion of Theorem~\ref{tk propositon}, we have
\begin{align*}
{F}_i(x^{k_j})&\preceq {F}_i(x^k),\ \forall i=\overline{1,m},\\
\Rightarrow\  {F}_i(\overline{x})&\preceq {F}_i(x^k)\ \text{as }j\to\infty,\ \forall i=\overline{1,m}.
\end{align*}
Thus, the sequence $\{{F}_i(x^k)\}$ is bounded below and consequently, by the second assertion of Theorem~\ref{tk propositon}, we have
\begin{align*}
\sum_{k=0}^{\infty}\|d^k\|^2<\infty,
\end{align*}
equivalently,
\[
\lim_{k\to\infty}d^k=0.
\]
In particular, for a subsequence $\{d^{k_j}\}$ of $\{d^k\}$, we have
\[
\lim_{j\to \infty}d^{k_j}=0.
\]
Further, since $d^{k_j}=\alpha(x^{k_j})$ and $\alpha(\cdot)$ is continuous on
$int(dom\ F)$ by the fifth assertion of Theorem~\ref{Theorem 1}, we obtain
\begin{align*}
\alpha(\overline{x})=\lim_{j\to\infty}\alpha(x^{k_j})=\lim_{j\to\infty}d^{k_j}=0 .
\end{align*}
Thus, it follows that
\begin{align*}
\beta(\overline{x})=\psi\big(\overline{x};\alpha(\overline{x})\big)+\tfrac12\|\alpha(\overline{x})\|^2
=\psi(\overline{x};0)=0.
\end{align*}
\noindent Hence, from the second assertion of 
Theorem~\ref{Theorem 1}, we conclude that $\overline{x}$ is a Pareto critical point of the
MIVOP$_C$. 
\end{proof}

\begin{corollary}
  Assume  that $dom\ F=\mathbb{R}^n$ and
the set $\mathcal{U}_0:=\{x\in dom\ F:\ F_i(x)\preceq F_i(x^0),\ i=\overline{1,m}\}$
is bounded. Then, the sequence $\{x^k\}$ generated by Algorithm~\ref{alg:1} is bounded, admits at least one accumulation
point and all of its accumulation points are Pareto critical points of MIVOP$_C$.
\end{corollary}
\begin{proof}
Since the sequence $\{F_i(x^k)\}$ is non-increasing, we have
\begin{align*}
F_i(x^{k+1})\preceq F_i(x^k),\ \forall i=\overline{1,m}\ \text{and}\ k\ge 0, 
\end{align*}
which upon induction and from Lemma~\ref{lemma ppt}, yields 
\[
F_i(x^k)\preceq F_i(x^0)\ \Rightarrow\ x^k\in\mathcal{U}_0,\ \forall k.
\]
Hence, the sequence $\{x^k\}$ is bounded and thus, by the Bolzano--Weierstrass theorem, it
has at least one accumulation point. Further, since $dom\ F=\mathbb{R}^n$, we have
$int(dom\ F)=\mathbb{R}^n$, therefore by Theorem~\ref{thm:global}, every accumulation point lies in $\R^n$ and each of them is a Pareto
critical point of the MIVOP$_C$.
\end{proof}

\noindent We now establish the global rate of convergence of Algorithm~\ref{alg:1} under the nonconvex setting.

\begin{theorem}
  Suppose that, for each $i\in\mathcal I\subseteq\{1,2,\dots,m\}$, the IVF $F_i(x^k)$ has a lower bound $F^{\min}_i\in \IR$ for all $k\ge 0$, and let $F^{\min}:=\displaystyle\min_{i\in\mathcal I}F_i^{\min}$ and $F_0^{\max}:=\displaystyle\max_{i=1,\ldots,m}{F_i}(x^0)$. Further, if $\{x^k\}$ is the sequence generated by Algorithm~\ref{alg:1}, then the following inequality holds:
   \[
   \min_{0\le j\le k-1}\|d^j\|\le \sqrt{\frac{\min \big\{\underline {F_0}^{\max}-\underline{F}^{\min},\overline{F_0}^{\max}-\overline{F}^{\min}\big \}}{\sigma\ t_{\min}}} \frac{1}{\sqrt{k}},
   \]
   where $d^j$ is the proximal gradient direction computed at the iterate $x^j$, $\sigma\in (0,1)$ and $t_{\min}$ is the uniform lower bound on the step length $t_k$ given in Theorem~\ref{t minimum}.
\end{theorem}

\begin{proof}
   Let $i\in\mathcal I\subseteq \{1,\dots,m\}$. From the Armijo-like condition~\eqref{armijo 2}, at the $k^{th}$-iteration we have
    \begin{equation*}
        \begin{aligned}
            \underline{F_i}(x^k+t_kd^k)&\le \underline{F_i}(x^k)+\sigma t_k \overline{h_i}(x^k;d^k),\\
            \overline{F_i}(x^k+t_kd^k)&\le \overline{F_i}(x^k)+\sigma t_k \overline{h_i}(x^k;d^k).
        \end{aligned}
    \end{equation*}
Since $\overline{h_i}(x^k;d^k)\le \psi(x^k;d^k)$, the above inequalities together with Lemma~\ref{Lemma 4.1}, yield
\begin{equation*}
   \begin{aligned}
     \underline{F_i}(x^k+t_kd^k)&\le \underline{F_i}(x^k)-\sigma t_k\|d^k\|^2 \\ 
     \text{and }\ \overline{F_i}(x^k+t_kd^k)&\le \overline{F_i}(x^k)-\sigma t_k \|d^k\|^2,\ \forall i\in \mathcal I
  \end{aligned} 
  \end{equation*} 
Further, since $t_k\ge t_{\min}$, it follows from above inequalities
\begin{equation}
\begin{aligned}
\underline{F_i}(x^k+t_kd^k)- \underline{F_i}(x^k)&\le- \sigma t_{\min} \|d^k\|^2\\
\text{and }\ \overline{F_i}(x^k+t_kd^k)-\overline{F_i}(x^k)&\le- \sigma t_{\min} \|d^k\|^2,\ \forall i\in\mathcal I.
\label{inequality}
\end{aligned}
\end{equation}
Now, on telescoping the first inequality from $j=0$ to $k-1$, we obtain
\begin{align}
    \sum_{j=0}^{k-1}\Big \{\underline{F_i}(x^{j+1})-\underline{F_i}(x^j) \Big \}&\le- \sigma\ t_{\min}\sum_{j=0}^{k-1}\|d^j\|^2, \notag \\
    \Rightarrow\  \underline{F_i}(x^k)-\underline{F_i}(x^0)&\le -(\sigma\  t_{\min})k\displaystyle\min_{0\le j\le k-1}\|d^j\|^2,\notag \\[4pt]
    \text{or, }\  k \displaystyle\min_{0\le j\le k-1}\|d^j\|^2&\le \frac{ \underline{F_i}(x^0)-\underline{F_i}(x^k)}{\sigma \ t_{\min}},\ \forall i\in \mathcal I.
    \label{rate lower Fi}
\end{align}
Similarly, the second inequality in \eqref{inequality} admits the same and therefore
\begin{align}
 k \displaystyle\min_{0\le j\le k-1}\|d^j\|^2\le \frac{ \overline{F_i}(x^0)-\overline{F_i}(x^k)}{\sigma \ t_{\min}},\ \forall i\in \mathcal I.
 \label{upper Fi rate}
\end{align}
Also, by the definitions of $F^{\min}$
and ${F_0}^{\max}$,
we have
\begin{align*}
    [\, \underline{F}^{\min},\overline{F}^{\min}\, ]&\preceq [\, \underline{F_i}(x^k),\overline{F_i}(x^k)\, ],\ \forall i\in \mathcal I\\
    \text{and }\ [\, \underline{F_i}(x^0),\overline{F_i}(x^0)\, ]&\preceq [\, \underline{F_0}^{\max},\overline{F_0}^{\max}\, ],\ \forall i=\overline{1,m}.
\end{align*}
Thus, it follows that
\begin{align*}
     \underline{F_i}(x^0)-\underline{F_i}(x^k)\le \underline{F_0}^{\max}-\underline{F}^{\min}\ \text{and}\  \overline{F_i}(x^0)-\overline{F_i}(x^k)\le \overline{F_0}^{\max}-\overline{F}^{\min},\ \forall i\in \mathcal I.
\end{align*}
Therefore, using the above inequalities in \eqref{rate lower Fi} and \eqref{upper Fi rate}, we obtain
\begin{align*}
k \displaystyle&\min_{0\le j\le k-1}\|d^j\|^2\le \frac{\underline{F_0}^{\max}-\underline{F}^{\min}}{\sigma \ t_{\min}}\ \text{ and }\ k\displaystyle\min_{0\le j\le k-1}\|d^j\|^2\le \frac{\overline{F_0}^{\max}-\overline{F}^{\min}}{\sigma \ t_{\min}}\\[4pt]
&\text{or, }\  \displaystyle\min_{0\le j\le k-1}\|d^j\|\le \sqrt{\frac{\min \big\{\underline {F_0}^{\max}-\underline{F}^{\min},\overline{F_0}^{\max}-\overline{F}^{\min}\big \}}{\sigma\ t_{\min}}} \frac{1}{\sqrt{k}}.
\end{align*}
 This completes the proof.
\end{proof}

\section{Numerical experiments}\label{Section 5}


In this section, we demonstrate the performance of Algorithm~\ref{alg:1} on 12 test problems, where in each problem, the smooth component $f_i\ (i=\overline{1,m})$ involved in the objective function is adopted from the corresponding test problems considered in Mondal and Ghosh \cite{mondal2025steepest} and to incorporate the nonsmooth component, we define the IVF $g_i\ (i=\overline{1,m})$ as an interval valued $\ell_1$-norm function, given as follows:
\[
g_i(x)= \frac{1}{(i+1)^n}\|x-i+1\|_1\odot [1,2],\ \forall i=1,2\dots,m,
\]
where $\|x\|_1=\sum_{j=1}^n|x_j|$, for every $x\in\R^n$. Note that for each test problem $dom\ F=\R^n$. 

\medskip 
\noindent All the computational results are obtained using MATLAB R2025b on a personal computer.  To execute Algorithm~\ref{alg:1} in MATLAB, the following initializations are made:


\begin{itemize}
    \item For the initial point, we use the MATLAB function `rand' that chooses random point within the bounds of the decision variable mentioned in Table~\ref{tb: test problem}.
    \item The Armijo-like line search strategy is implemented with fixed parameters $\gamma=0.5$ and $\sigma =0.001$ across all the test instances. 
    \item For obtaining the search direction $d^k$ and the optimal value $\beta^k$, the subproblem COP is solved using the MATLAB optimization toolbox `quadprog'. 
    \item From the steps of Algorithm~\ref{alg:1}, the stopping criteria is set as $\beta^k\ge-\epsilon$, where $\epsilon=10^{-6}$.
\end{itemize}

\noindent The following Table~\ref{tb: test problem} summarizes the 12 test problems which reports the dimension of the decision space ($n$), the number of objectives ($m$), the componentwise lower and the upper limits of the decision variables ($Lb$ and $Ub$), and the nature of the corresponding problem. The notation I-XYZ$_{\ell_1}$ indicates the test problem I-XYZ given in Mondal and Ghosh~\cite{mondal2025steepest} together with an interval valued $\ell_1$-norm function described as above.

\begin{table}[H]
\centering
\footnotesize
\label{tab:testproblems}
\begin{tabular}{cccccc}
\toprule
Test Problem & $n$ & $m$ &  $Lb^\top$ & $Ub^\top$ & Nature \\

\midrule

I-BK1$_{\ell_1}$ & 2 & 2 & $(-5,\,-5)$  & $(10,\,10)$ & Convex \\

I-VU2$_{\ell_1}$ & 2 & 2 & $(-4,\,-4)$  & $(4,\,4)$ & Convex \\

I-CH$_{\ell_1}$  & 2 & 2 & $(-5,\,-5)$  & $(5,\,5)$ & Convex \\

I-KW2$_{\ell_1}$ & 2 & 2 & $(-3,\,-1)$  & $(0,\,2)$ & Non-convex \\

I-Far1$_{\ell_1}$ & 2 & 2 & $(-1,\,-1)$  & $(1,\,1)$ & Non-convex \\

I-SD$_{\ell_1}$  & 4 & 2 & $(1,\,\sqrt{2},\,\sqrt{2},\,1)$  & $(6,\,6,\,6,\,6)$ & Convex \\

I-IKK1$_{\ell_1}$ & 2 & 3 & $(-50,\,-50)$  & $(50,\,50)$ & Convex \\

I-AP1$_{\ell_1}$ & 2 & 3 & $(-100,\,-100)$  & $(100,\,100)$ & Convex \\

I-MOP7$_{\ell_1}$ & 2 & 3 & $(-400,\,-400)$  & $(400,\,400)$ & Convex \\

I-MHHM2$_{\ell_1}$ & 2 & 3 & $(0,\,0)$  & $(1,\,1)$ & Convex \\

I-Comet$_{\ell_1}$ & 3 & 3 & $(1,\,-2,\,0)$  & $(3.5,\,2,\,1)$ &  Non-convex \\

I-VFM2$_{\ell_1}$ & 3 & 3 & $(-5,\,-5,\,-5)$  & $(10,\,10,\,10)$ & Convex \\

\bottomrule
\end{tabular}
\caption{Summary of the test problems.}
\label{tb: test problem}
\end{table}

\noindent To illustrate the implementation of Algorithm~\ref{alg:1}, 
we consider the I-BK1$_{\ell_1}$ test problem listed in Table~\ref{tb: test problem}. Accordingly, the associated MIVOP$_C$ is given as follows:
\begin{align*}
    \textbf{MIVOP1}_C\qquad \min_{x\in\R^2} F(x):=(F_1(x),\ F_2(x))^\top,
\end{align*}
where each of the component $F_i(x)=f_i(x)\oplus g_i(x)$, is given as
\begin{align*}
    F_1(x)&=[0.1x_1^2+0.1x_2^2,\ 0.2x_1^2+0.3x_2^2]\oplus \dfrac{\|x\|_1}{4}\odot [1,2],\\
    F_2(x)&=[0.1(x_1-5)^2+0.1(x_2-5)^2,\ 0.3(x_1-5)^2+0.5(x_2-5)^2]\oplus \dfrac{\|x-1\|_1}{9}\odot[1,2].
\end{align*}
We choose an initial point $x^0=(8.569794,\,-3.706053)^\top$ from the decision boudary $Lb^\top\le (x_1,x_2)^\top\le Ub^\top$ and using the initializations mentioned above, the results obtained are presented in Table~\ref{tab:IBK1_{OP}1}.



\begin{table}[H]
\centering
\footnotesize
\renewcommand{\arraystretch}{1.15}
\setlength{\tabcolsep}{2.5pt}
\begin{tabular}{c c c c c}
\toprule
$k$ & $x^k=(x_1^k,x_2^k)^\top$ & $F_1(x^k)$ & $F_2(x^k)$ & $\beta(x^k)$ \\
\midrule
0  & $(8.569794,\,-3.706053)$ & $[11.78658,\;24.94665]$ & $[10.21786,\;44.44868]$ & $-1.950566\times10^{0}$ \\
1  & $(7.378578,\,-2.130574)$ & $[8.27556,\;17.00506]$ & $[6.70684,\;29.23298]$ & $-1.248362\times10^{0}$ \\
2  & $(6.425605,\,-0.870190)$ & $[6.02851,\;12.13275]$ & $[4.45979,\;19.46056]$ & $-7.861105\times10^{-1}$ \\
3  & $(5.717369,\,0.042077)$ & $[4.70887,\;9.41792]$ & $[3.14015,\;13.70606]$ & $-4.038637\times10^{-1}$ \\
4  & $(5.024623,\,0.614641)$ & $[3.97228,\;7.98233]$ & $[2.41320,\;10.59586]$ & $-2.218211\times10^{-1}$ \\
5  & $(4.481198,\,0.999777)$ & $[3.47831,\;7.05658]$ & $[2.01392,\;8.85529]$ & $-4.594205\times10^{-2}$ \\
6  & $(4.253138,\,1.199291)$ & $[3.31586,\;6.77554]$ & $[1.88392,\;8.15724]$ & $-1.602683\times10^{-2}$ \\
7  & $(4.119279,\,1.318184)$ & $[3.22997,\;6.63371]$ & $[1.81508,\;7.77447]$ & $-5.618971\times10^{-3}$ \\
8  & $(4.040286,\,1.388881)$ & $[3.18258,\;6.55806]$ & $[1.77714,\;7.55844]$ & $-1.974207\times10^{-3}$ \\
9  & $(3.993550,\,1.430882)$ & $[3.15569,\;6.51613]$ & $[1.75565,\;7.43417]$ & $-6.944356\times10^{-4}$ \\
10 & $(3.965860,\,1.455825)$ & $[3.14017,\;6.49228]$ & $[1.74325,\;7.36180]$ & $-2.444322\times10^{-4}$ \\
11 & $(3.949442,\,1.470634)$ & $[3.13110,\;6.47849]$ & $[1.73602,\;7.31933]$ & $-8.606589\times10^{-5}$ \\
12 & $(3.939702,\,1.479426)$ & $[3.12578,\;6.47043]$ & $[1.73177,\;7.29430]$ & $-3.030711\times10^{-5}$ \\
13 & $(3.933924,\,1.484645)$ & $[3.12264,\;6.46569]$ & $[1.72927,\;7.27950]$ & $-1.067678\times10^{-5}$ \\
14 & $(3.930495,\,1.487742)$ & $[3.12078,\;6.46289]$ & $[1.72778,\;7.27074]$ & $-3.759702\times10^{-6}$ \\
15 & $(3.928460,\,1.489581)$ & $[3.11968,\;6.46124]$ & $[1.72691,\;7.26554]$ & $-1.323347\times10^{-6}$ \\
16 & $(3.927252,\,1.490672)$ & $[3.11902,\;6.46026]$ & $[1.72639,\;7.26247]$ & $-4.665929\times10^{-7}$ \\
\bottomrule
\end{tabular}
\caption{Iterations generated by Algorithm~\ref{alg:1} for the I-BK1$_{\ell_1}$ problem.}
\label{tab:IBK1_{OP}1}
\end{table}

\noindent It is  observed that the proposed method generates a sequence of iterates that converges to an approximate Pareto critical point $\hat x=(3.927252,1.490672)^\top$ in 16 iterations. 
Further, since the chosen test problem I-BK1$_{\ell_1}$ is convex, therefore, from Proposition~\ref{ppsn 1}, the point $\hat x=(3.927252,1.490672)^\top$ is a weakly Pareto optimal point of MIVOP1$_C$.

\begin{figure}[H]
\centering
\includegraphics[width=0.5\linewidth]{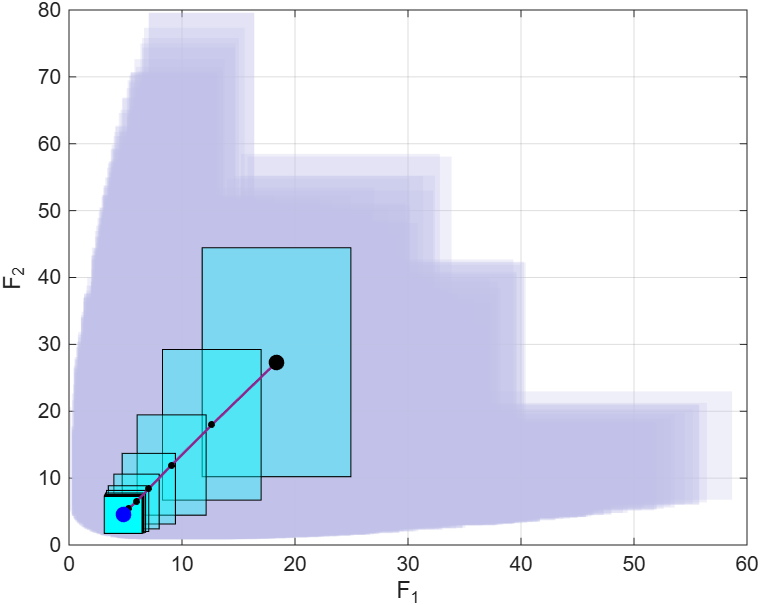}
\caption{Graphical representation of Algorithm~\ref{alg:1} on I-BK1$_{\ell_1}$ problem.}
\label{fig: IBK1}
\end{figure}

\medskip
\noindent The graphical representation of the trajectory of the sequence $\{F(x^k)\}$ generated by Algorithm~\ref{alg:1} is depicted in  Figure~\ref{fig: IBK1}. The gray region depicts the feasible region of the objective function of the I-BK1$_{\ell_1}$ problem that is approximated by the union of objective values, $i.e$, $\displaystyle\bigcup_{Lb\le x\le Ub}F(x)$, corresponding to 1500 uniformly distributed random sample points generated over the decision space.
The cyan rectangles denote the objective values $F(x^k)$ generated by the iterates of Algorithm~\ref{alg:1}, while the purple line joining their centers illustrates the trajectory of the iterates in the objective feasible space. The trajectory starts from the black filled circles that represent the center of the rectangle $F(x^k)$ and end to a blue filled circle that represent the center of the rectangle $F(\hat x)$, where $\hat x$ is the Pareto critical point of the MIVOP1$_C$.

\begin{figure}[H]
\centering
\begin{subfigure}{0.32\textwidth}
    \centering
    \includegraphics[width=\linewidth]{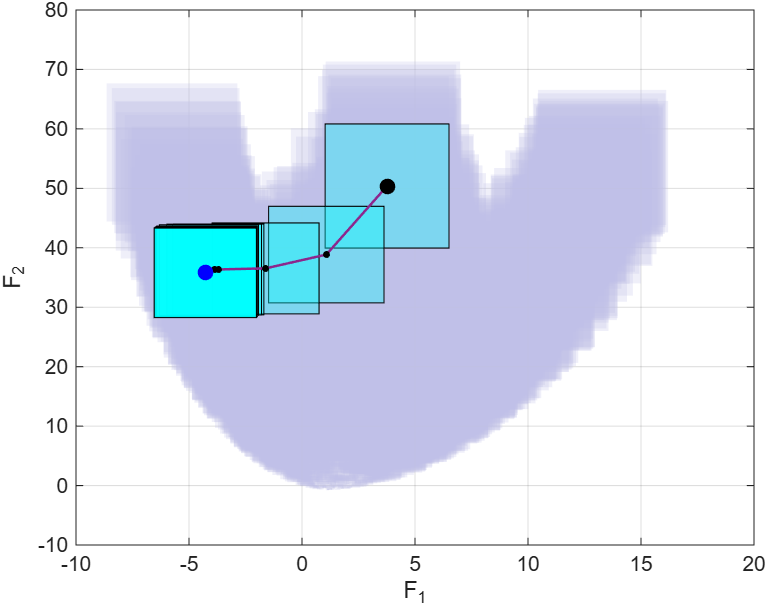}
    \caption{I-VU2$_{\ell_1}$}
\end{subfigure}
\hfill
\begin{subfigure}{0.32\textwidth}
    \centering
    \includegraphics[width=\linewidth]{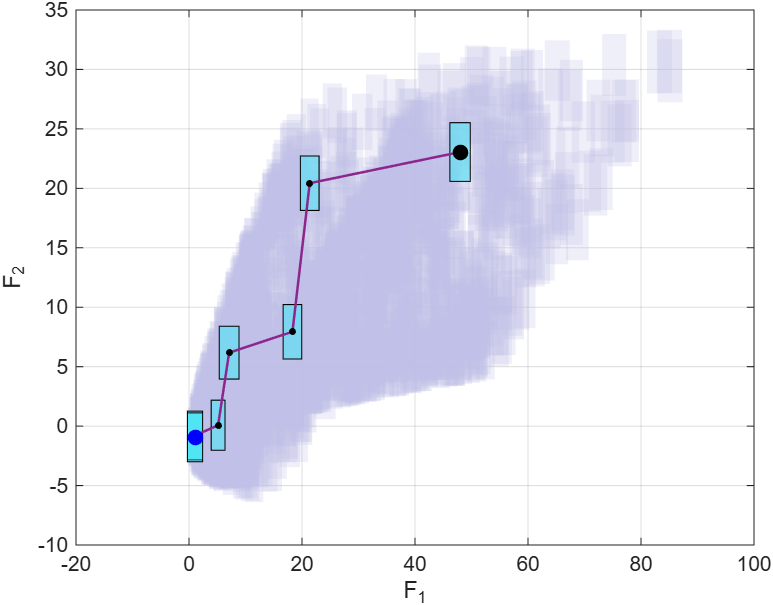}
    \caption{I-CH$_{\ell_1}$}
\end{subfigure}
\hfill
\vspace{0.5cm}
\begin{subfigure}{0.32\textwidth}
    \centering
    \includegraphics[width=\linewidth]{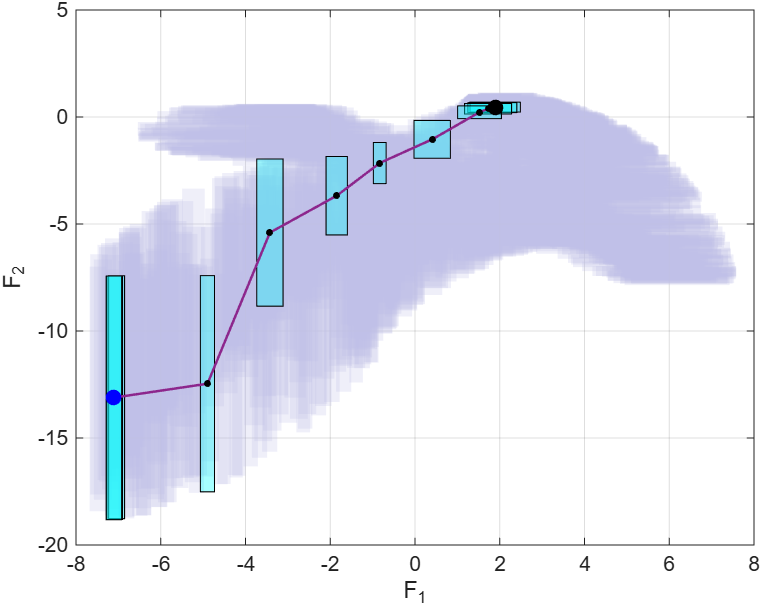}
    \caption{I-KW2$_{\ell_1}$}
\end{subfigure}
\hfill

\vspace{0.5cm}

\begin{subfigure}{0.32\textwidth}
    \centering
    \includegraphics[width=\linewidth]{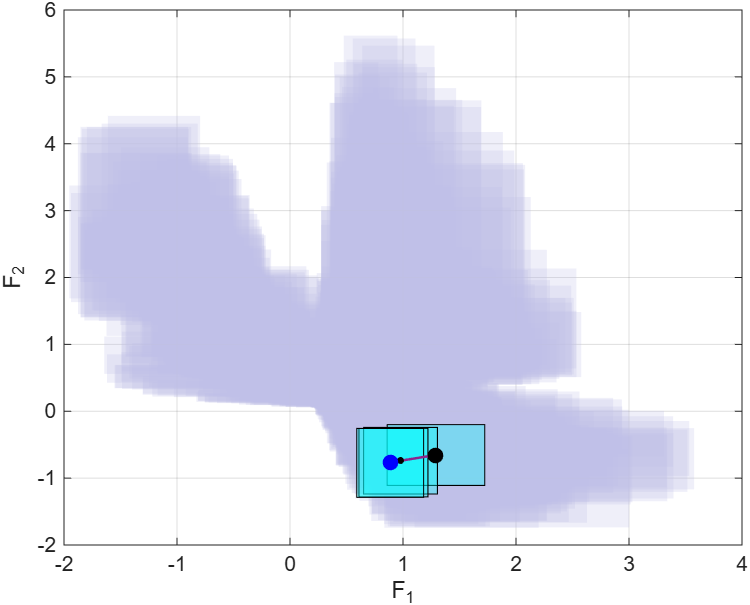}
    \caption{I-Far1$_{\ell_1}$}
\end{subfigure}
\hspace{1.5cm}
\begin{subfigure}{0.32\textwidth}
    \centering
   \includegraphics[width=\linewidth]{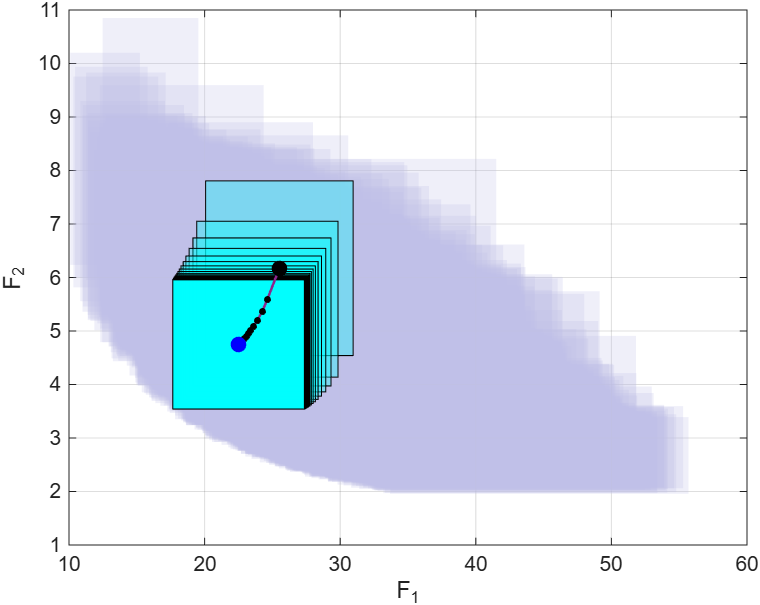}
    \caption{I-SD$_{\ell_1}$}
\end{subfigure}

\caption{Trajectory of biobjective $F(x^k)$ generated by Algorithm~\ref{alg:1}.}
\label{fig: trajectory}
\end{figure}

\noindent In a similar manner, the graphical representation of Algorithm~\ref{alg:1} on other five bi-objective ($m=2$) test problems listed in Table~\ref{tb: test problem}, starting from a random initial point to a Pareto critical point, is depicted in Figure~\ref{fig: trajectory}.

\begin{figure}[H]
\centering
\begin{subfigure}{0.32\textwidth}
    \centering
    \includegraphics[width=\linewidth]{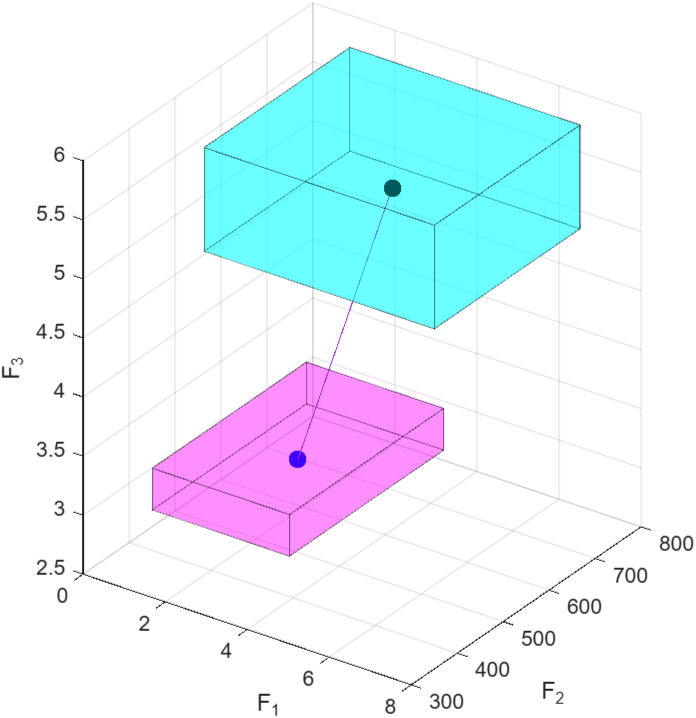}
    \caption{I-IKK1$_{\ell_1}$}
\end{subfigure}
\hfill
\begin{subfigure}{0.32\textwidth}
    \centering
    \includegraphics[width=\linewidth]{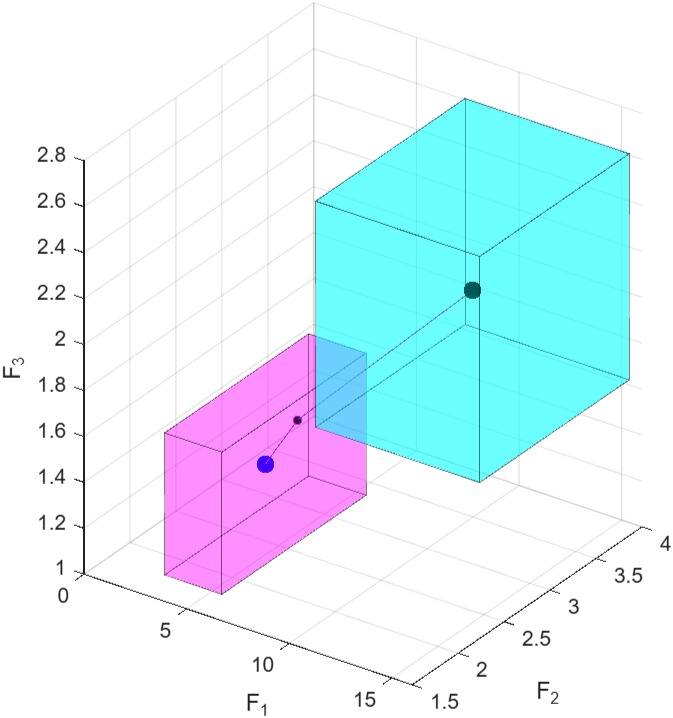}
    \caption{I-AP1$_{\ell_1}$}
\end{subfigure}
\hfill
\begin{subfigure}{0.32\textwidth}
    \centering
    \includegraphics[width=\linewidth]{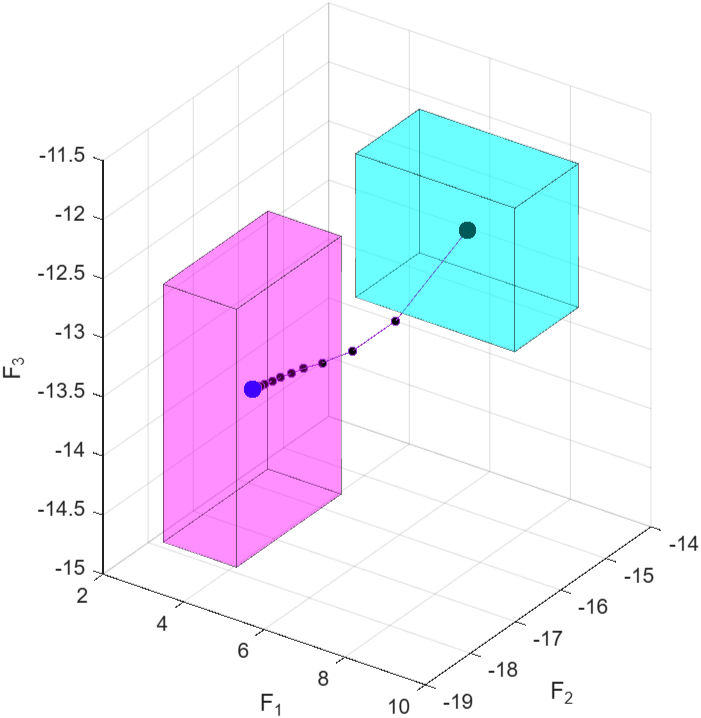}
    \caption{I-MOP7$_{\ell_1}$}
\end{subfigure}
\hfill

\vspace{0.5cm}

\begin{subfigure}{0.32\textwidth}
    \centering
    \includegraphics[width=\linewidth]{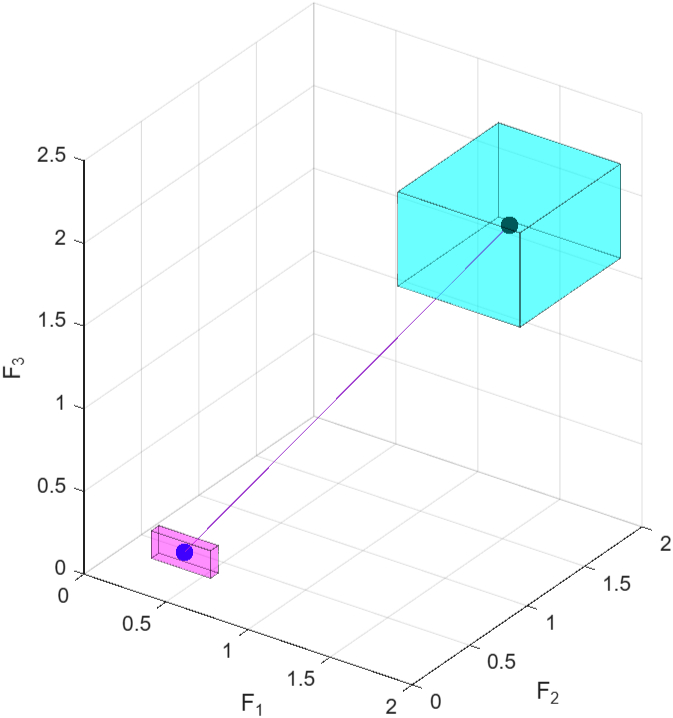}
    \caption{I-MHHM2$_{\ell_1}$}
\end{subfigure}
\hfill
\begin{subfigure}{0.32\textwidth}
    \centering
    \includegraphics[width=\linewidth]{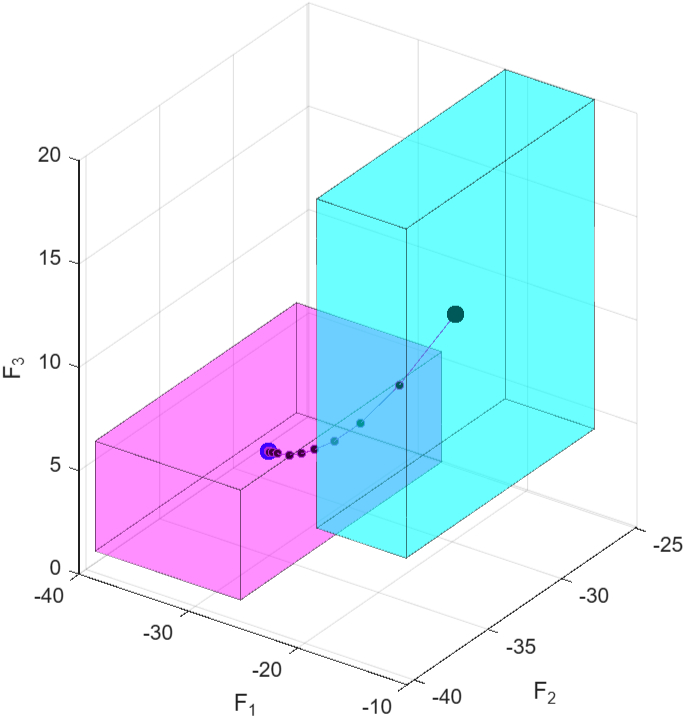}
    \caption{I-Comet$_{\ell_1}$}
\end{subfigure}
\hfill
\begin{subfigure}{0.32\textwidth}
    \centering
    \includegraphics[width=\linewidth]{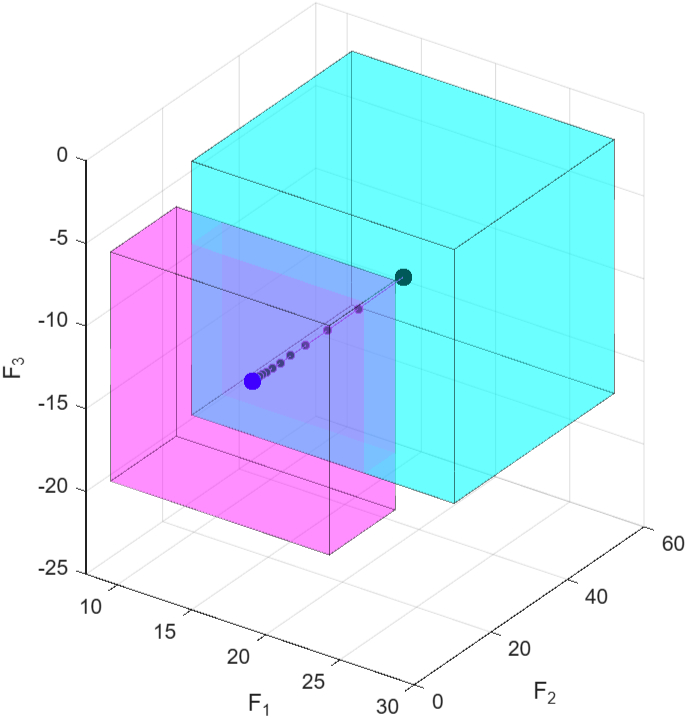}
    \caption{I-VFM2$_{\ell_1}$}
\end{subfigure}

\caption{Trajectory of tri-objective $F(x^k)$ generated by Algorithm~\ref{alg:1}.}
\label{fig: trajectory tri}
\end{figure}

\medskip
\noindent In case of the tri-objective ($m=3$) test problems listed in Table~\ref{tb: test problem}, visualizing the complete objective feasible region is a difficult task. Hence, for an initial point $x^0$ chosen randomly, we have displayed only the objective values $F(x^0)$ and $F(\hat x)$ in Figure~\ref{fig: trajectory tri}. The cyan cuboid denotes the objective value at the initial point, while the magenta cuboid corresponds to the objective value at the Pareto critical point. The purple curve, extending from the black point to the blue point, represents the trajectory of the sequence $\{F(x^k)\}$. 

\noindent Now, to assess the computational efficiency and consistency of Algorithm~\ref{alg:1}, we conduct a numerical  evaluation using two performance metrics: iteration count and the CPU execution time. The proposed method is run for 100 random initial points for all the test problems listed in Table~\ref{tb: test problem} and the maximum (max), minimum (min), mean ($\mu$), median (Med), and standard deviation (SD) of both the metrics are evaluated. The results of the performance of Algorithm~\ref{alg:1} based on these two metrics are given in Table~\ref{tab: performance}.

\begin{table}[H]
\centering
\footnotesize
\renewcommand{\arraystretch}{1.15}
\setlength{\tabcolsep}{4.5pt}
\begin{tabular}{ccc}
\toprule

Test Problem  & 

\shortstack{Iterations \\ (Max, Min, $\mu$, Med, SD)} & 

\shortstack{CPU time(s) \\ (Max, Min, $\mu$, Med, SD) } \\

\midrule

I-BK1$_{\ell_1}$   & $(16,\,0,\,9.7200,\,11.0000,\,4.9565)$ & $(1.2494,\,0.0011,\,0.0260,\,0.0153,\,0.1237)$  \\

I-VU2$_{\ell_1}$   & $(30,\,0,\,8.1300,\,5.5000,\,7.9209)$ & $(0.0347,\,0.0010,\,0.0105,\,0.0076,\,0.0090)$  \\

I-CH$_{\ell_1}$     & $(12,\,0,\,6.8600,\,7.0000,\,2.5147)$ & $(0.0134,\,0.0020,\,0.0081,\,0.0083,\,0.0090)$  \\

I-SD$_{\ell_1}$   & $(85,\,0,\,33.4700,\,33.5000,\,16.1235)$ & $(0.1384,\,0.0011,\,0.0452,\,0.0448,\,0.0223)$  \\

I-KW2$_{\ell_1}$   & $(19,\,0,\,5.1200,\,5.000,\,4.2977)$ & $(0.0270,\,0.0009,\,0.0076,\,0.0067,\,0.0057)$ \\

I-Far1$_{\ell_1}$  & $(29,\,0,\,2.0900,\,1.0000,\,4.0703)$ & $(0.0373,\,0.0009,\,0.0036,\,0.0021,\,0.0050)$ \\

I-IKK1$_{\ell_1}$ & $(299,\,0,\,88.9200,\,42.0000,\,99.8384)$ &
$(0.4012,\,0.0011,\,0.1182,\,0.0560,\,0.1317)$ \\

IAP1$_{\ell_1}$ & $(299,\,0,\,127.4600,\,53.500,\,136.1506)$ &
$(0.5527,\,0.0020,\,0.1979,\,0.0877,\,0.2104)$ \\

I-MOP7$_{\ell_1}$ & $(299,\,24,\,208.8300,\,299.0000,\,121.0983)$ & 
$(1.1415,\,0.0388,\,0.6763,\,0.9948,\,0.4577)$ \\

I-MHHM2$_{\ell_1}$ & $(3,\,0,\,1.6200,\,2.0000,\,0.5646)$ & 
$(0.0068,\,0.011,\,0.0034,\,0.0036,\,0.0009)$ \\

I-VFM2$_{\ell_1}$ & $(299,\,0,\,124.4600,\,17.0000,\,140.3510)$ & 
$(1.2593,\,0.0013,\,0.4347,\,0.0296,\,0.5211)$ \\

I-Comet$_{\ell_1}$ & $(165,\,0,\,9.4600,\,7.5000,\,17.0782)$ &
$(0.2500,\,0.0012,\,0.0166,\,0.0147,\,0.0259)$ \\

\bottomrule
\end{tabular}
\caption{Performance of Algorithm~\ref{alg:1}.}
\label{tab: performance}
\end{table}

\noindent \textbf{Note: } To the best of our knowledge, no existing numerical algorithm has been specifically developed for the MIVOP$_C$ considered in this work. Therefore, instead of benchmarking against unavailable methods, we have evaluated the proposed algorithm through a collection of numerical test problems that demonstrate its convergence and practical performance.

\section{Application}\label{section 6}

In this section, we illustrate the applicability of Algorithm~\ref{alg:1} on a mean-variance bi-objective portfolio model with transaction cost. In the mean--variance framework of Markowitz~\cite{markowitz1952}, an investor seeks to simultaneously maximize the expected portfolio return and minimize the portfolio risk, measured through the variance of returns. However,  in practice, rebalancing an existing portfolio incurs a transaction cost that is typically modeled through a nonsmooth, piecewise-linear term in the decision variables, resulting in a composite objective structure. 

\medskip
\noindent Let us consider a portfolio containing two assets, and let $x_1$ and $x_2$ be the fractions of total investment funds allocated to the first and second asset, respectively. 
Let $\mathcal R_1$ and $\mathcal R_2$ denote the expected returns, and $Cov(\mathcal R_1,\mathcal R_1)$ and  $Cov(\mathcal R_2,\mathcal R_2) $ represent the  variances of the returns  $\mathcal R_1$ and $\mathcal R_2$ of the two assets, respectively. Further,   $Cov(\mathcal R_1,\mathcal R_2)$ represents the covariances between the returns such that $Cov(\mathcal R_1,\mathcal R_2)=Cov(\mathcal R_2,\mathcal R_1)$. 

\medskip
\noindent In addition to the standard return and risk components, we also look for the transaction cost incurred when the investor rebalance an existing portfolio. 
Let $\overline{x}_1$ and $\overline{x}_2$ denote the current holding of the two assets, respectively, and $\mathcal C_1$ and $\mathcal C_2$ be the transaction cost rate, charged per unit of the two assets bought or sold. 

\medskip
\noindent The main goal of the investor has two folds: first to maximize the expected portfolio return net of transaction cost, and second to minimize the portfolio risk. This leads to the following multiobjective portfolio optimization problem:
\begin{align*}
    \textbf{MPO}_{tc}\quad  &\min F(x_1,x_2)=(F_1(x_1,x_2),F_2(x_1,x_2))^\top\\
    &\text{subject to }x_1+x_2=1,\\
    &\qquad \quad\quad \ \  x_1,x_2\ge 0.
\end{align*}
where each component $F_i(x_1,x_2)=f_i(x_1,x_2)+g_i(x_1,x_2)$, for $i=1,2$, are defined as follows:
\begin{align*}
    f_1(x_1,x_2)&=-(\mathcal R_1 x_1+\mathcal R_2x_2),\ \ g_1(x_1,x_2)=\mathcal C_1|x_1-\overline{x}_1|+\mathcal C_2|x_2-\overline{x}_2|,\\
    f_2(x_1,x_2)&= Cov(\mathcal R_1,\mathcal R_1)x_1^2+2Cov(\mathcal R_1,\mathcal R_2)x_1 x_2+Cov(\mathcal R_2,\mathcal R_2)x_2^2\ \text{ and }\   g_2(x_1,x_2)=0.
\end{align*}

\noindent In practice, the expected returns and the associated portfolio risk (covariance of asset returns) cannot be known with certainty, as they fluctuate due to market volatility, estimation error, and incomplete information. To account for this uncertainty, we represent these parameters using closed intervals, where the lower and upper endpoints specify the range of possible values.

\medskip
\noindent For instance, we take $\mathcal R_1=[4,5],\ \mathcal R_2=[2,3],\ Cov(\mathcal R_1,\mathcal R_1)=[1,2],\ Cov(\mathcal R_1,\mathcal R_2)=[0.1,0.3],\ Cov(\mathcal R_2,\mathcal R_2)=[1.5,3]$. Further, due to uncertainty in trading costs caused by bid-ask spreads and market impacts, we assume $\mathcal C_1=[0.008,0.012]$ and $\mathcal C_2=[0.006,0.009]$ and the investor's current holdings are taken as $(\overline{x}_1,\overline{x}_2)^\top=(0.6,0.4)^\top$. Based on these parameter, the above MPO$_{tc}$ problem in terms of multiobjective interval valued optimization problem is given as follows:
\begin{align*}
    \textbf{MIVOP2}_{C}\qquad &\min F(x_1, x_2):=\big (F_1(x_1,x_2), F_2(x_1,x_2)\big )^\top\\
    &\text{subject to }x_1+x_2=1,\\
    &\qquad \qquad \ \  x_1,x_2\ge 0,
\end{align*}
where each of the component $F_i(x_1,x_2)=f_i(x_1,x_2)\oplus g_i(x_1,x_2)$, for $i=1,2,$ are
\begin{align*}
    F_1(x_1,x_2)&=\big((-1)\odot \big([4,5]\odot x_1\oplus[2,3]\odot x_2\big )\big )\oplus\big ([0.008,0.012]\odot |x_1-0.6|\\
   &\qquad  \hspace{5.9cm} \oplus [0.006,0.009]\odot |x_2-0.4|\big ),\\
    F_2(x_1,x_2)&=\big ([1,2]\odot x_1^2\oplus[0.1,0.3]\odot 2x_1x_2\oplus[1.5,3]\odot x_2^2\big )\oplus \mathbf{0}.
\end{align*}
Since $x_1+x_2=1$, therefore the above problem can be transformed to a single variable unconstrained optimization problem, as follows:
\begin{align*}
  \textbf{MIVOP3}_{C}\ & \min_{x_1\in \R_+}\ F(x_1,x_2):= \Big(\big( (-1)\odot [2x_1+2,2x_1+3]\big )\oplus\big ([0.014,0.021]\odot|x_1-0.6|\big ),\\
    &\qquad \qquad \hspace{3.3cm}[2.3 x_1^2-2.8x_1+1.5,4.4x_1^2-5.4x_1+3]\Big )
\end{align*}

\noindent Applying Algorithm~\ref{alg:1} on MIVOP3$_{C}$ with different initial values $x_1^0\in \{0,0.25,0.50,0.75,1.0\}$ and considering the initializations of the parameters same as taken in Section~\ref{Section 5}, the  results obtained are presented in the following Table~\ref{tab: MIPO solution}. 

\begin{table}[H]
\centering
\footnotesize
\renewcommand{\arraystretch}{1.10}
\setlength{\tabcolsep}{2.3pt}
\begin{tabular}{cccc}
\toprule
$x_1^0$ & Pareto critical point $(\hat x_1,\hat x_2)^\top$ & $F_1(\hat x)$ & $F_2(\hat x)$   \\

\midrule

$0$ & $(0.507,\,0.493)$ & $[-4.012698,\,-3.012047]$ & $[ 0.671613,\, 1.393216]$   \\

$0.25$ & $(0.675,\,0.325)$ & $[-4.348950,\,-3.348425]$ & $[ 0.657938,\, 1.359750]$ \\

$0.50$  & $(0.50,\,0.50)$ & $[-3.998600,\,-2.997900]$ & $[ 0.675000,\, 1.400000]$ \\

$0.75$  & $(0.75,\,0.25)$ & $[-4.497900,\,-3.496850]$ & $[ 0.693750,\, 1.425000]$ \\

$1$ & $(1,\,0)$ & $[-4.994400,\,-3.991600]$ & $[ 1.000000,\, 2.000000]$ \\

\bottomrule
\end{tabular}
\caption{Results of Algorithm~\ref{alg:1} for MIVOP3$_C$.}
\label{tab: MIPO solution}
\end{table}

\section{Conclusion}\label{section 7}

In this article, we proposed a proximal gradient method for unconstrained multiobjective interval-valued optimization problems with composite objective functions (MIVOP$_ C$). We introduced an auxiliary IVF to characterize the Pareto criticality of the MIVOP$_ C$ and based on these characterizations, we formulated a strongly convex subproblem whose unique solution determines the proximal gradient direction. We, then developed an implementable algorithm with an Armijo-type line search that is compatible with the interval partial order. Under suitable  assumptions, we derived the existence of a uniform lower bound on the step lengths and established that every sequence generated by the proposed method converges globally to a Pareto critical point of MIVOP$_C$ together with a global convergence rate of order 
$\mathcal{O}(1/\sqrt{k})$ in the nonconvex setting. Further, the numerical experiments on test problems demonstrated the effectiveness of the method, while its application to a mean-variance portfolio selection problem with transaction costs, illustrated its practical applicability to optimization problems involving uncertainty and nonsmooth terms.

\medskip
\noindent For future research, we will focus on establishing improved convergence rates  under suitable convexity and strong convexity assumptions. Furthermore, it would be worthwhile to develop accelerated or proximal Newton-type methods for the considered problem and to extend the proposed framework to constrained problems.\\

\noindent \textbf{Author  Contributions} The authors contributed equally to  this study.\\

\noindent \textbf{Funding} No funding was received for conducting this study.\\

\noindent \textbf{Data Availability} No datasets were generated or analyzed in this study.

\section*{Declarations}

\textbf{Competing interests} The authors declare no competing interests.

\bibliographystyle{plain}

\end{document}